\documentclass[a4paper,10pt]{article}
\usepackage{stmaryrd}
\usepackage{amsfonts}
\usepackage{bbm}
\usepackage{amscd}
\usepackage{mathrsfs}
\usepackage{latexsym,amssymb,amsmath,amscd,amscd,amsthm,amsxtra,xypic}
\usepackage[dvips]{graphicx}

\newtheorem{thm}{Theorem}[section]
\newtheorem{lem}[thm]{Lemma}
\newtheorem{cor}[thm]{Corollary}
\newtheorem{pro}[thm]{Proposition}
\newtheorem{ex}[thm]{Example}
\newtheorem{rmk}[thm]{Remark}
\newtheorem{defi}[thm]{Definition}

\newcommand{\lon }{\,\rightarrow\,}
\newcommand{\be }{\begin{eqnarray*}}
\newcommand{\ee }{\end{eqnarray*}}

\newcommand{\defbe}{\triangleq}
\newcommand{\pf}{\noindent{\bf Proof.}\ }

\newcommand{\Real}{\mathbb R}

\newcommand{\huaE}{\mathcal{E}}
\newcommand{\huaF}{\mathcal{F}}
\newcommand{\huaG}{\mathcal{G}}

\newcommand{\huaK}{\mathcal{K}}

\newcommand{\CWM}{C^{\infty}(M)}

\newcommand{\set}[1]{\left\{#1\right\}}

\newcommand{\frkd}{\mathfrak d}

\newcommand{\frkg}{\mathfrak g}

\newcommand{\frkr}{\mathfrak r}

\newcommand{\frkt}{\mathfrak t}

\newcommand{\frky}{\mathfrak y}

\newcommand{\frkL}{\mathfrak L}

\newcommand{\frkX}{\mathfrak X}

\def\gpd{\,\lower1pt\hbox{$\longrightarrow$}\hskip-.24in\raise2pt
         \hbox{$\longrightarrow$}\,}

\def\qed{\hfill ~\vrule height6pt width6pt depth0pt}

\newcommand{\rhowx}{\rho^{\star}}

\newcommand{\Liederivative}{\frkL}

\newcommand{\half}{\frac{1}{2}}

\newcommand{\pairing}[1]{\left\langle #1\right\rangle}
\newcommand{\ppairingE}[1]{\left ( #1\right )_E}
\newcommand{\starpair}[1]{\left ( #1\right )_{*}}

\newcommand{\pibracket}[1]{\left [ #1\right ]_{\pi}}
\newcommand{\conpairing}[1]{\left\langle  #1\right\rangle }
\newcommand{\Courant}[1]{\left\llbracket  #1\right\rrbracket }
\newcommand{\Dorfman}[1]{\{ #1\}}
\newcommand{\CDorfman}[1]{   [    #1   ]_{\huaK}   }
\newcommand{\jetc}[1]{   [    #1   ]_{\jet{C}}   }

\newcommand{\jet}{\mathfrak{J}}

\newcommand{\jetd}{\mathbbm{d}}
\newcommand{\dev}{\mathfrak{D}}

\newcommand{\TStarM}{T^*M}

\newcommand{\pie}{^\prime}
\newcommand{\Id}{\mathbf{1}}

\newcommand{\e}{\mathbbm{e}}
\newcommand{\p}{\mathbbm{p}}
\newcommand{\id}{\mathbbm{i}}
\newcommand{\jd}{\mathbbm{j}}

\newcommand{\dM}{\mathrm{d}}
\newcommand{\dA}{\mathrm{d}^A}
\newcommand{\dB}{\mathrm{d}^B}
\newcommand{\dimension}{\mathrm{dim}}
\newcommand{\omni}{\mathcal{E}}

\newcommand{\Hom}{\mathrm{Hom}}
\newcommand{\Der}{\mathrm{Der}}
\newcommand{\Ad}{\mathrm{Ad}}
\newcommand{\Aut}{\mathrm{Aut}}
\newcommand{\gl}{\mathfrak {gl}}

\newcommand{\Img}{\mathrm{Im}}

\newcommand{\rhoAWuXing}{\rho_A^{\star}}
\newcommand{\rhoBWuXing}{\rho_B^{\star}}
\begin{document}
\title{
{ $E$-Courant Algebroids
\thanks
 {
 Research partially supported by NSFC grants 10871007 and 10911120391/
A0109. The third author is also supported by CPSF grant 20090451267.
 }
} }
\author{ Zhuo Chen$^1$, ~ Zhangju Liu$^2$ and Yunhe Sheng$^3$\\
$^1$Department of Mathematics, \\Tsinghua University, Beijing 100084, China\\
$^{2}$Department of Mathematics and LMAM\\ Peking University,
Beijing 100871, China\\
$^{3}$Department of Mathematics \\ Jilin University,
 Changchun 130012, Jilin, China\\
          {\sf email: $^1$zchen@math.tsinghua.edu.cn,\quad
          $^2$liuzj@pku.edu.cn,\quad $^3$shengyh@jlu.edu.cn} }

\date{}


\maketitle
\begin{abstract}
In this paper, we introduce the notion of  $E$-Courant algebroids,
where $E$ is a vector bundle. It is a   kind of generalized Courant
algebroid  and contains   Courant algebroids, Courant-Jacobi
algebroids and omni-Lie algebroids as its special cases.  We explore
novel phenomena exhibited by $E$-Courant algebroids and   provide
many examples. We study the automorphism groups of  omni-Lie
algebroids and  classify  the isomorphism classes of exact
$E$-Courant algebroids. In addition, we introduce the concepts  of
$E$-Lie bialgebroids and Manin triples.
\end{abstract}
\tableofcontents

\section{Introduction}
In recent years,
 Courant algebroids are  widely studied  from several aspects.
 They are applied in many mathematical
 objects such as
   Manin pairs and moment
maps \cite{AK,Courant morphism,Xu1,Libland and M}, generalized
complex structures \cite{BCG,GualtieriGeneralizedComplex,Xu2},
$L_{\infty}$-algebras and symplectic supermanifolds \cite{ROy},
gerbes \cite{PW}, BV algebras and  topological field theories
\cite{BV quantization,ROy1}.

We recall two notions closely related to Courant algebroids ---
Jacobi bialgebroids and omni-Lie algebroids. Jacobi bialgebroids and
generalized Lie bialgebroids are introduced, respectively,  in
\cite{Grabowski marmo1} and \cite{generalize Lie bialgebroid}   to
generalize Dirac structures from Poisson manifolds to Jacobi
manifolds. More general geometric objects are
  generalized Courant algebroids \cite{generalized}  and Courant-Jacobi algebroids
 \cite{Grabowski marmo2}.
The notion of omni-Lie algebroids, a generalization of the notion of
omni-Lie algebras introduced in \cite{Weinomni},   is defined in
\cite{clomni} in order to characterize all possible Lie algebroid
structures on a vector bundle $E$.  An omni-Lie algebra can be
regarded as the linearization of the exact Courant algebroid
$TM\oplus \TStarM$ at a point and is studied from several aspects
recently \cite{BCG,InteCourant,shengliuzhuomni-Lie 2,UchinoOmni}.
Moreover, Dirac structures of omni-Lie algebroids are studied by the
authors in \cite{CLS}.

In this paper, we introduce a kind of generalized Courant algebroid
 called  $E$-Courant
algebroids. The values of the anchor map of an $E$-Courant algebroid
lie in $\dev E$, the  bundle of differential operators. Moreover,
its Dirac structures are necessarily Lie algebroids equipped with a
representation on $E$. The notion of $E$-Courant algebroids  not
only unifies Courant-Jacobi algebroids and omni-Lie algebroids, but
also provides a number of interesting objects, e.g. the
$T^*M$-Courant algebroid structure on the jet bundle of a Courant
algebroid over $M$ (Theorem \ref{thm:jet C}).

Recall that an exact Courant algebroid structure  on $TM\oplus T^*M$
is a
  twist of the standard Courant algebroid  by a closed 3-form \cite{PW}.
  This structure  includes twisted
Poisson structures  and is related to gerbes and  topological sigma
models  \cite{twistP,Poisson quasi-Nijenhuis}. In this paper, we are
inspired to study exact $E$-Courant algebroids
 similar to the situation of exact Courant algebroids.

 We also study the automorphism groups   of  omni-Lie algebroids, for which
 we need the language of Leibniz cohomologies \cite{Loday and
Pirashvili,Leibniz diff mfd}. Moreover, we   introduce the notion of
$E$-Lie bialgebroids, which  generalizes the notion of generalized
Lie bialgebroids. We shall prove that, for an $E$-Lie bialgebroid,
there induces on the underlying vector bundle $E$  a Lie algebroid
structure $(\mathrm{rank}(E)\geq 2),$ or a local Lie algebra
structure $(\mathrm{rank}(E)=1)$ (Theorem \ref{Thm:ElocalLie}).

This  paper is organized  as follows. In Section
\ref{Sec:GeneralizedCourant}  we introduce  the notion of
$E$-Courant algebroids.  We prove that the jet bundle $\jet C$ of a
Courant algebroid $C$ over $M$ admits a natural $T^*M$-Courant
algebroid structure.
 In Section 3 we discuss  the properties of
  $E$-dual pairs of Lie algebroids.
   In Section 4  we find  the automorphism
groups  and all possible twists  of  omni-Lie algebroids. In Section
5 we  study exact $E$-Courant algebroids and   prove that every
exact $E$-Courant algebroid with an isotropic splitting is
isomorphic to an omni-Lie algebroid. In general, an exact
$E$-Courant algebroid is a twist of the standard omni-Lie algebroid
by a 2-cocycle in the Leibniz cohomology of $\Gamma(\dev E)$ with
coefficients in $\Gamma(\jet E)$,  which can also be treated as a
3-cocycle in the Leibniz cohomology of $\Gamma(\dev E)$ with
coefficients in $\Gamma(E)$. In Section 6 we  study $E$-Lie
bialgebroids. In Section 7 we extend the theory of Manin triples
from the context of Lie bialgebroids to $E$-Lie bialgebroids and
give some interesting examples.

 {\bf Acknowledgement:} ~~ Z. Chen  would like to  give his warmest
thanks to  P. Xu and M. Grutzmann for their useful comments.  Y.
Sheng gives his warmest thanks to L. Hoevenaars, M. Crainic, I.
Moerdijk and C. Zhu for their useful comments during his stay in
Utrecht University  and Courant Research Center, G\"{o}ttingen. We
also give our warmest thanks to the referees for many helpful
suggestions and pointing out typos and erroneous statements.


\section{$E$-Courant
algebroids}\label{Sec:GeneralizedCourant}

 Let
$E{\lon}M$ be a vector bundle and $\dev E$ the associated covariant
differential operator bundle. Known as the gauge Lie algebroid of
the
 frame bundle
 $\huaF(E)$ (see \cite[Example
3.3.4]{Mkz:GTGA}),  $\dev E$ is a transitive Lie algebroid with the
Lie bracket $[\cdot,\cdot]_\dev$ (commutator).  The corresponding
Atiyah sequence is as follows:
\begin{equation}\label{Seq:DE}
\xymatrix@C=0.5cm{0 \ar[r] & \gl(E)  \ar[rr]^{\id} &&
                \dev{E}  \ar[rr]^{\jd} && TM \ar[r]  & 0.
                }
\end{equation}
In \cite{clomni}, the authors proved that  the jet bundle $\jet E$
(see \cite{CFSecondary,Jet bundle} for more details about jet
bundles) can be regarded as an $E$-dual bundle of $\dev E$, i.e.
\begin{eqnarray*}  {\jet E}  &\cong &
\set{\nu\in \Hom( \dev{E}  ,E )\,|\,  \nu(\Phi)=\Phi\circ
\nu(\Id_E),\quad\forall ~~ \Phi\in \gl(E )}  \subset \Hom( \dev{E}
,E ).
\end{eqnarray*}
Associated to  the jet bundle $\jet E$,  the jet sequence of $E$ is
given by:
\begin{equation}\label{Seq:JetE}
\xymatrix@C=0.5cm{0 \ar[r] & \Hom(TM,E) \ar[rr]^{\quad\quad\quad\e}
&&
                {\jet}{E} \ar[rr]^{\p} && E \ar[r]  & 0.
                }
\end{equation}
The operator $\jetd: \Gamma(E) \rightarrow \Gamma(\jet E)$ is given
by:
$$\jetd u (\frkd) := \frkd (u),\quad \quad \forall~  u  \in   \Gamma(E),  ~~~\frkd \in\Gamma(\dev E).$$
The following formula  is needed.
\begin{equation}\label{eqn:d fX} \jetd(fX)=\dM f\otimes X+f\jetd
X,\quad \forall~X\in\Gamma(C),~f\in C^\infty(M).
\end{equation}
For a vector bundle $\huaK$ over $M$ and a bundle map
$\rho:\huaK\longrightarrow \dev E$, we denote the induced
$E$-adjoint bundle map by $\rho^\star$, i.e.
\begin{equation}\label{E adjoint map}
\rho^\star : \Hom(\dev E,E)\rightarrow \Hom(\huaK,E), \quad
\rho^\star(\nu)(k)=\nu(\rho(k)), \quad \forall~k\in
\huaK,~\nu\in\Hom(\dev E,E).
\end{equation}
The notion of Leibniz algebras is introduced by Loday
\cite{Loday,Loday and Pirashvili,Leibniz foliation}. A Leibniz
algebra $\frkg$ is an $R$-module, where $R$ is a commutative ring,
endowed with a linear map
$[\cdot,\cdot]:\frkg\otimes\frkg\longrightarrow\frkg$ satisfying
$$
[g_1,[g_2,g_3]]=[[g_1,g_2],g_3]+[g_2,[g_1,g_3]],\quad
\forall~g_1,g_2,g_3\in \frkg.
$$

\begin{defi}\label{defi: $E$-C A}
An $E$-Courant algebroid is a quadruple $(\huaK,
\ppairingE{\cdot,\cdot},\CDorfman{\cdot,\cdot},\rho)$, where
\begin{enumerate}\item[$\bullet$]
$\huaK$ is a vector bundle  over $M$ such that
$(\Gamma(\huaK),\CDorfman{\cdot,\cdot})$ is a Leibniz
algebra;\item[$\bullet$] $\ppairingE{\cdot,\cdot}:
~\huaK\otimes\huaK\lon E$ is a symmetric nondegenerate  $E$-valued
pairing, which induces an embedding: $\huaK \hookrightarrow
\Hom(\huaK,E)$; \item[$\bullet$] the anchor $\rho:~ \huaK\lon \dev
E$ is a bundle map, \end{enumerate}such that the following
properties hold for all $~X,Y,Z\in \Gamma(\huaK)$:
\begin{itemize}
\item[~~~~ \rm (EC-1)]  \quad \quad \quad \quad
$\rho\CDorfman{X,Y} =[\rho(X),\rho(Y)]_{\dev}$;
\item[~~~~ \rm (EC-2)] \quad \quad \quad \quad  $\CDorfman{X,X}  =  \rhowx\jetd\ppairingE{X,X}
$;
\item[~~~~ \rm (EC-3)]\quad \quad \quad \quad
$\rho(X)\ppairingE{Y,Z}=\ppairingE{\CDorfman{X,Y} ,Z}
+\ppairingE{Y,\CDorfman{X,Z} } $;
 \item[~~~~ \rm (EC-4)]  \quad \quad \quad \quad $\rhowx( \jet E)\subset \huaK$,
i.e. $
 \ppairingE{\rhowx(\mu),~X}=\half\mu(\rho(X)),~ ~~\forall~ \mu\in\jet
 E;
 $
\item[~~~~ \rm (EC-5)] \quad \quad \quad \quad
$\rho\circ\rhowx=0$.
\end{itemize}
\end{defi}
\begin{rmk}
 If the $E$-valued pairing
$\ppairingE{\cdot,\cdot}:\huaK\otimes\huaK\longrightarrow E$ is
surjective,  Properties (EC-4) and (EC-5) can be inferred from
Property (EC-2). In particular, if $E$ is a line bundle, any
nondegenerate $E$-valued pairing $\ppairingE{\cdot,\cdot}$ is
surjective.
\end{rmk}
\begin{lem}For any $X,Y\in\Gamma(\huaK)$ and $f\in C^\infty(M)$, we
have
\begin{eqnarray}
\label{eqn:bracket X fY}\CDorfman{X, fY}&=&f\CDorfman{X,
Y}+(\jd\circ\rho(X)f)Y,\\
\label{eqn:bracket fX Y}\CDorfman{fX, Y}&=&f\CDorfman{X,
Y}-(\jd\circ\rho(Y)f)X+2\rhowx(\dM f\otimes\ppairingE{X,Y}).
\end{eqnarray}
\end{lem}
\pf By Property (EC-3), for all $X,~Y,~Z\in \Gamma(\huaK)$ and
$~f\in\CWM$, we have
\begin{eqnarray*}
\ppairingE{\CDorfman{X,fY} ,Z}
+\ppairingE{fY,\CDorfman{X,Z} }&=&\rho(X)\ppairingE{fY,Z}\\
&=&\jd\circ\rho(X)(f)\ppairingE{Y,Z}+f\rho(X)\ppairingE{Y,Z}\\
&=&\jd\circ\rho(X)(f)\ppairingE{Y,Z}+f\ppairingE{\CDorfman{X,Y},Z}
+f\ppairingE{Y,\CDorfman{X,Z}}.
\end{eqnarray*}
Since the pairing $\ppairingE{\cdot,\cdot}$ is nondegenerate, it
follows
 that
$$
\CDorfman{X,fY}=\jd\circ\rho(X)(f)Y+f\CDorfman{X,Y}.
$$
By Property (EC-2), we have
$$
\CDorfman{X, fY}+\CDorfman{fY,
X}=2\rhowx\jetd(f\ppairingE{X,Y})=2f\rhowx\jetd\ppairingE{X,Y}+2\rhowx(\dM
f\otimes\ppairingE{X,Y}).
$$
Substitute $\CDorfman{X, fY}$ by (\ref{eqn:bracket X fY}) and apply
Property (EC-2) again, we obtain (\ref{eqn:bracket fX Y}). \qed
\vspace{3mm}

For a subbundle $L\subset\huaK$, denote by $L^\bot\subset\huaK$ the
subbundle
$$
L^\bot=\{e\in\huaK\mid~\ppairingE{e,l}=0,~\forall~l\in L\}.
$$
\begin{defi}
A Dirac structure of an $E$-Courant algebroid $(\huaK,
\ppairingE{\cdot,\cdot},\CDorfman{\cdot,\cdot},\rho)$ is a subbundle
$L\subset \huaK$ which is   closed under the bracket
$\CDorfman{\cdot,\cdot}$ and satisfies $L=L^\bot$.
\end{defi}

Evidently, $L=L^\bot$ implies that $L$ is maximal isotropic with
respect to the $E$-valued pairing $\ppairingE{\cdot,\cdot}$. In
general,   $L$ being maximal isotropic with respect to
$\ppairingE{\cdot,\cdot}$ does not imply
 $L=L^\bot$.
\begin{ex}
Let $\huaK$ be $\mathbb R^3$  with the standard  basis
$e_1,~e_2,~e_3$. The $\mathbb R$-valued pairing
$(\cdot,\cdot)_{\mathbb R}$ is given by
$$
(e_1,e_3)_{\mathbb R}=(e_2,e_2)_{\mathbb R}=1,\quad
(e_1,e_1)_{\mathbb R}=(e_1,e_2)_{\mathbb R}=(e_2,e_3)_{\mathbb
R}=(e_3,e_3)_{\mathbb R}=0.
$$
Obviously, $L=\mathbb Re_1$ is maximal isotropic but $L^\bot=\mathbb
Re_1\oplus\mathbb Re_2 \neq L$.
\end{ex}

\begin{pro}\label{pro:rep}
Any Dirac structure $L$ has an induced   Lie algebroid structure and
is equipped with a Lie algebroid representation
$\rho_L=\rho\mid_L:~L\lon\dev E$ on $E$.
\end{pro}
\pf Given a Dirac structure $L$, by Property (EC-2), we  have
$\CDorfman{X,X}=0,$ for all $~ X\in\Gamma(L)$, which implies that
$\CDorfman{\cdot,\cdot}\mid_L$ is skew-symmetric. By
(\ref{eqn:bracket X fY}),
$(L,\CDorfman{\cdot,\cdot}\mid_{L},(\jd\circ\rho)\mid_L)$ is a Lie
algebroid. Finally by Property (EC-1), $\rho_L:~L\lon\dev E$ is
  a representation.  \qed
\begin{rmk}
If $E$ is the trivial line bundle $M\times\mathbb R$, then $\dev
E\cong TM\oplus (M\times \mathbb R)$. Thus we can decompose
$\rho=a+\theta$, for some $a:\huaK\longrightarrow TM$ and
$\theta:\huaK\longrightarrow
 M\times\mathbb R$. For a Dirac structure $L$, since $\rho_L $ is a representation of the Lie algebroid
$L$,  it follows that $\theta_L=\theta|_L \in \Gamma(L^*)$ is a
1-cocycle in the Lie algebroid cohomology of $L$. Therefore,
$(L,\theta_L)$ is a Jacobi algebroid, which is, by definition, a Lie
algebroid $A$ together with a 1-cocycle $\theta\in\Gamma(A^*)$ in
the Lie algebroid cohomology \cite{Grabowski marmo2}.
\end{rmk}

 One may refer to \cite{Mkz:GTGA} for more general theories of Lie
algebroids,  Lie algebroid cohomologies and their representations.
Now we briefly recall the notions of  omni-Lie algebroids,
generalized Courant algebroids,
 Courant-Jacobi algebroids and
  generalized Lie bialgebroids. We will see that  $E$-Courant
  algebroids
 unify all these structures.

$\bullet$ \textbf{Omni-Lie algebroids}

The notion of omni-Lie algebroids is introduced in \cite{clomni} to
characterize Lie algebroid structures on a vector bundle. It is a
generalization of Weinstein's omni-Lie algebras. Recall that there
is a natural symmetric nondegenerate $E$-valued pairing
$\conpairing{\cdot,\cdot}_E$ between $\jet E$ and $\dev{E}$:
\begin{eqnarray}\nonumber
\conpairing{\mu,\frkd}_E=\conpairing{\frkd,\mu}_E &\defbe& \frkd
u,\quad\forall ~~ \mu=[u]_m\in {\jet
E},~u\in\Gamma(E),~\frkd\in\dev{E}.
\end{eqnarray}
Moreover, this pairing  is $C^\infty(M) $-linear and satisfies the
following properties:
\begin{eqnarray*}
\conpairing{\mu, \Phi }_E &=& \Phi\circ \p(\mu),\quad\forall ~ ~\Phi\in \gl(E),~\mu\in{\jet E};\\
\conpairing{ {\frky} ,\frkd}_{E} &=& {\frky}\circ
\jd(\frkd),\quad\forall ~~ \frky\in \Hom(TM,E),~\frkd\in\dev{E}.
\end{eqnarray*}
Furthermore, $\Gamma (\jet E)$ is  invariant  under any Lie
derivative   $\Liederivative_{\frkd}$,
 $\frkd \in\Gamma(\dev{E})$, which is defined by the
 Leibniz rule:
\begin{equation}\label{def:Lie d}
\conpairing{\Liederivative_{\frkd}\mu,\frkd\pie}_{E}\defbe
\frkd\conpairing{\mu,\frkd\pie}_{E}-\conpairing{\mu,[\frkd,\frkd\pie]_{\dev}}_{E},
\quad\forall~ \mu \in \Gamma(\jet{E}), ~
~\frkd\pie\in\Gamma(\dev{E}).
\end{equation}

\begin{defi}\label{def:omni algebroid}{\rm\cite{clomni}}
Given a vector bundle $E$, the quadruple
$(\omni,\Dorfman{\cdot,\cdot},\ppairingE{\cdot,\cdot},\rho)$ is
called the {\em omni-Lie algebroid} associated to  $E$, where $
\omni=\dev{E}\oplus \jet{E}$, the anchor $\rho$ is the projection
from $\omni$ to $\dev{E}$, the bracket operation
$\Dorfman{\cdot,\cdot}$ and the nondegenerate $E$-valued pairing
$\ppairingE{\cdot,\cdot}$ are given respectively by
\begin{eqnarray}
\label{standard pair}\ppairingE{\frkd+\mu,\frkr+\nu}&\defbe&
\half(\conpairing{\frkd,\nu}_E
+\conpairing{\frkr,\mu}_E),\\
\label{standard bracket}\Dorfman{\frkd+\mu,\frkr+\nu}&\defbe&
[\frkd,\frkr]_{\dev}+\Liederivative_{\frkd}\nu-\Liederivative_{\frkr}\mu
+ \jetd\conpairing{\mu,\frkr}_E\,.
\end{eqnarray}
\end{defi}
If there is no risk of confusion, we simply denote the omni-Lie
algebroid
$(\omni,\Dorfman{\cdot,\cdot},\ppairingE{\cdot,\cdot},\rho)$ by
$\omni$.
 We call the $E$-valued pairing (\ref{standard pair}) and
the bracket   (\ref{standard bracket}), respectively, {\bf{ the
standard pairing }} and {\bf {the standard bracket}} on $\omni=\dev
E\oplus\jet E$. One may refer to \cite{clomni} for more details of
the property of  omni-Lie algebroids. Evidently, the $E$-adjoint map
$\rhowx$ is $\Id_{\jet E}$, the identity map on $\jet E$. It is
easily seen that the omni-Lie algebroid $\omni$ is an $E$-Courant
algebroid. Its Dirac structures are  studied  by the authors in
\cite{CLS}. \vspace{3mm}

$\bullet$ \textbf{Generalized Courant algebroids (Courant-Jacobi
algebroids)}

The notion of generalized Courant algebroids is   introduced in
\cite{generalized}. It is a pair $(\huaK,\rho^\theta)$ subject to
some compatibility conditions, where $\huaK\lon M$ is a vector
bundle
 equipped with a nondegenerate symmetric bilinear form
$(\cdot,\cdot)$, a skew-symmetric bracket $[\cdot,\cdot]$ on
$\Gamma(\huaK)$ and a bundle map $\rho^\theta:\huaK\lon TM\times
\mathbb R$, which is also a first-order differential operator. We
may write $\rho^\theta(X)=(\rho(X),\langle\theta,X\rangle)$, where
$\rho:~\huaK\lon TM$ is linear and $\theta\in\Gamma(\huaK^*)$
satisfies
$$
\theta([X,Y])=\rho(X)\theta(Y)-\rho(Y)\theta(X),\quad
\forall~X,Y\in\Gamma(\huaK).
$$
 One should note that the skew-symmetric
bracket $[\cdot,\cdot]$ does not satisfy the Jacobi identity.  The
notion of Courant-Jacobi algebroids is introduced in \cite{Grabowski
marmo2}. In \cite{generalized}, it is established the equivalence of
generalized Courant algebroids  and
 Courant-Jacobi algebroids. Roughly speaking, the difference between
them is that the    generalized Courant algebroid has a
skew-symmetric bracket $[\cdot,\cdot]$ and a Courant-Jacobi
algebroid has an operation $\circ$, which is also known as the
Dorfman bracket \cite{Dorfman1993}. The former does not satisfy the
Jacobi identity, while the later satisfies the Leibniz rule.
Moreover, $[\cdot,\cdot]$ can be realized as the skew-symmetrization
of $\circ$.  A generalized Courant algebroid reduces to a Courant
algebroid
 if $\theta=0$ (see \cite{LWXmani}).

Evidently, a generalized Courant algebroid  is  an $E$-Courant
algebroid if we take $E=M\times \mathbb R$. It follows that all
Jacobi algebroids and Courant algebroids are  $M\times \mathbb
R$-Courant algebroids.\vspace{3mm}

$\bullet$ \textbf{Generalized Lie bialgebroids}

A Lie bialgebroid is a pair of vector bundles in duality, each of
which is a Lie algebroid, such that the differential defined by one
of them on the exterior algebra of its dual is a derivative of the
Schouten bracket  \cite{Y,MackenzieX:1994}. A generalized Lie
bialgebroid \cite{generalize Lie bialgebroid}, or a Jacobi
bialgebroid \cite{Grabowski marmo1}, is a pair
$((A,\phi_0),(A^*,X_0))$, where $A$ and $A^*$ are two vector bundles
in duality,  and, respectively, equipped with  Lie algebroid
structures $(A,[\cdot,\cdot],a)$ and $(A^*,[\cdot,\cdot]_*,a_*)$.
The data $\phi_0\in\Gamma(A^*)$ and $X_0\in\Gamma(A)$ are 1-cocycles
in their respective Lie algebroid cohomologies such that for all
$X,Y\in\Gamma(A)$, the following conditions are satisfied:
\begin{eqnarray}
\label{eqn:gene Lie bi 1}d_{*X_0}[X,Y]=[d_{*{X_0}}X,Y]_{\phi_0}+[X,d_{*{X_0}}Y]_{\phi_0},\\
\label{eqn:gene Lie bi 2}\phi_0(X_0)=0,\quad
a(X_0)=-a_*(\phi_0),\quad\Liederivative_{*\phi_0}X+\Liederivative_{X_0}
X=0,
\end{eqnarray}
where  $d_{*X_0}$  is the $X_0$-differential of $A$,
$[\cdot,\cdot]_{\phi_0}$  is the $\phi_0$-Schouten bracket,
$\Liederivative_*$ and $\Liederivative$ are  the usual Lie
derivatives. For more information of these notations, please refer
to \cite{generalize Lie bialgebroid}. For a Jacobi manifold
$(M,X,\Lambda)$, $((TM\times \mathbb R,(0,1)),(T^*M\times \mathbb
R,(-X,0)))$ is a generalized Lie bialgebroid. Furthermore, for a
generalized Lie bialgebroid, there is an induced Jacobi structure on
the base manifold $M$. In particular, both $((A,\phi_0)$ and
$(A^*,X_0))$ are Jacobi algebroids. If $\phi_0=0$ and $X_0$=0, a
generalized Lie bialgebroid reduces to a Lie bialgebroid. It is
known that for a generalized Lie bialgebroid
$((A,\phi_0),(A^*,X_0)),$ there is a natural generalized Courant
algebroid   $(A\oplus A^*,\phi_0+X_0)$. \vspace{3mm}

We give more examples of $E$-Courant algebroids.

\begin{ex} \label{ex: A A dual tensor E}\rm {Let $A$ be a Lie algebroid
and $\rho_A:A\lon \dev E$ a representation of $A$ on a vector bundle
$E$. Let $\huaK=A\oplus (A^*\otimes E)$.  For any $X,~Y\in
\Gamma(A)$, $~\xi\otimes u,~ \eta\otimes v\in\Gamma(A^*\otimes E)$,
we define the following operations:
\begin{eqnarray*}
\rho(X+\xi\otimes u)&=& \rho_A(X), \\
\CDorfman{X+\xi\otimes u,Y+\eta\otimes
v}&=&[X,Y]+\Liederivative_{X}(\eta\otimes
v)-\Liederivative_{Y}(\xi\otimes u)+\rho_A^\star\circ\jetd(\langle Y,\xi\rangle u), \\
\ppairingE{X+\xi\otimes u,Y+\eta\otimes v}&=&\half(\langle
X,\eta\rangle v+\langle Y,\xi\rangle u).
\end{eqnarray*}
Evidently, $\rhowx=\rho_A^\star:\jet E\longrightarrow A^*\otimes E$
and it is straightforward to check that $(A\oplus (A^*\otimes E),
\CDorfman{\cdot,\cdot},\ppairingE{\cdot,\cdot},\rho)$ is an
$E$-Courant algebroid.} In \cite{Libland},  the notion of
$AV$-Courant algebroids is introduced in order to study generalized
CR structures, which is closely related to this example but twisted
by a 3-cocycle in the cohomology of  the Lie algebroid
representation $\rho_A$.
\end{ex}
\begin{ex}\rm{
Consider an $E$-Courant algebroid $\huaK$ whose anchor $\rho$ is
zero. Thus   $\rhowx=0$, and the bracket $\CDorfman{\cdot,\cdot}$ is
skew-symmetric. So $\huaK$ is a
 bundle of Lie algebras.  Property (EC-3) shows that there is an
 invariant $E$-valued pairing.
 We conclude that an $E$-Courant algebroid $\huaK$ whose anchor $\rho$ is
zero is equivalently a
 bundle of Lie algebras with an invariant $E$-valued pairing.}
\end{ex}
\begin{ex}\label{ex:g+V}\rm{An omni-Lie algebra $\gl(V)\oplus V$ is a
 special  omni-Lie algebroid  whose base manifold is a point,
 hence a  $V$-Courant algebroid. Moreover, one may consider
a Lie algebra $(\frkg,[\cdot,\cdot]_\frkg)$ with faithful
representation $\rho_\frkg:\frkg\longrightarrow\gl(V)$ on a vector
space $V$. This representation is called {\em nondegenerate} if for
any $v\in V $, there is some $A\in\frkg$ such that
$\rho_\frkg(A)(v)\neq 0$.  Introduce a nondegenerate $V$-valued
pairing $(\cdot,\cdot)_V$ and a bilinear bracket $[\cdot,\cdot]$ on
the space $\frkg\oplus V$:
\begin{eqnarray*}
(A+u,B+v)_V&=&\frac{1}{2}(\rho_\frkg(A)(v)+\rho_\frkg(B)(u)),\\
{[A+u,B+v]}&=&[A,B]_\frkg+\rho_\frkg(A)(v),\quad
\forall~A+u,~B+v\in\frkg\oplus V,
\end{eqnarray*}
where $\rho:\frkg\oplus V\longrightarrow\gl(V)$ is defined by
$\rho(A+u)=\rho_\frkg(A)$ for   $A+u\in \frkg\oplus V$. Following
from
\begin{equation}\label{eqn:ex rho wx}
\rhowx(u)(B+v)=\half\rho_\frkg(B)(u)=(u,B)_V,
\end{equation}
we have $\rhowx=\Id_V$, as a map $\jet V=V\longrightarrow V$.
Clearly, $(\frkg\oplus V,(\cdot,\cdot)_V,[\cdot,\cdot],\rho)$ is a
$V$-Courant algebroid.

The bracket defined above   appeared  in \cite{InteCourant}. For any
representation $\rho:\frkg\longrightarrow\gl(V)$, we call
$(\frkg\oplus V,[\cdot,\cdot])$   a {\bf{hemisemidirect}} product of
$\frkg$ with $V$. There is also a natural exact Courant algebra
 associated to  any $\frkg$-module \cite{BCG}. }
\end{ex}

The above  example can be generalized to the situation of Lie
algebroids.
\begin{ex}{\rm
Let $(A,[\cdot,\cdot],a)$ be a Lie algebroid with a nondegenerate
representation $\rho_A:A\longrightarrow \dev E$. On  the vector
bundle $A\oplus\jet E$, define an $E$-valued pairing
$\ppairingE{\cdot,\cdot}$ and a bracket $\{\cdot,\cdot\}$ by
\begin{eqnarray*}
\ppairingE{X+\mu,Y+\nu}&=&\half(\conpairing{\rho_A(X),\nu}_E+\conpairing{\rho_A(Y),\mu}_E),\\
\{X+\mu,Y+\nu\}&=&[X,Y]+\Liederivative_{\rho(X)}\nu-\Liederivative_{\rho(Y)}\mu+\jetd\conpairing{\rho_A(Y),\mu}_E,
\end{eqnarray*}
for any $X+\mu,~Y+\nu\in \Gamma(A\oplus\jet E)$, and define
$\rho:A\oplus\jet E\longrightarrow\dev E$ by
$\rho(X+\mu)=\rho_A(X)$. Similar
 to (\ref{eqn:ex rho wx}), we have $\rhowx=\Id_{\jet E}$. Then, it is
easily seen   that $(A\oplus\jet
E,\ppairingE{\cdot,\cdot},\{\cdot,\cdot\},\rho)$ is an $E$-Courant
algebroid.}
\end{ex}
 $\bullet$ \textbf{ The jet bundle of a Courant algebroid}

At the end of this section, we prove that for any Courant algebroid
$C$, $\jet C$ is a $T^*M$-Courant algebroid.
 The original
definition of a Courant algebroid is introduced in \cite{LWXmani}.
Here we use the alternative definition given by D. Roytenberg in
\cite{ROy}, that a Courant algebroid is a vector bundle
$C\longrightarrow M$ together with some compatible structures --- a
nondegenerate bilinear form $\pairing{\cdot,\cdot}$ on the bundle, a
bilinear operation $\Courant{\cdot,\cdot}$ on $\Gamma(E)$ and a
bundle map $a:~C\longrightarrow TM$ satisfying $a\circ a^*=0$. In
particular,  $(\Gamma(C),\Courant{\cdot,\cdot})$ is a Leibniz
algebra.

On the jet bundle $\jet C$ of the vector bundle $C$,  we introduce
the $T^*M$-valued pairing $\starpair{\cdot,\cdot}$, the bracket
$\jetc{\cdot,\cdot}$ and the anchor
$\rho:\jet{C}\longrightarrow\dev{(T^*M)}$ as follows.
\begin{enumerate}\item[a)]
  For any $X,~Y\in\Gamma(C)$,  the $T^*M$-valued pairing
$\starpair{\cdot,\cdot}$ of $\jetd X,~\jetd Y$ is given by
\begin{equation}\label{eqn:pair}
\starpair{\jetd X,\jetd Y}=\dM\pairing{X,Y}.
\end{equation}
By (\ref{eqn:d fX}),  we get
\begin{eqnarray*}
\starpair{\jetd X,\dM f\otimes Y}&=&\pairing{X,Y}\dM f,\\
\starpair{\dM f\otimes X,\dM f\otimes Y}&=&0.
\end{eqnarray*}

\item[b)] For any $X,~Y\in\Gamma(C)$,  the bracket $\jetc{\cdot,\cdot}$
of $\jetd X,~\jetd Y$ is given by
\begin{equation}\label{eqn:bracket}
\jetc{\jetd X,\jetd Y}=\jetd\Courant{X,Y}.
\end{equation}
By  (\ref{eqn:bracket X fY}), (\ref{eqn:bracket fX Y}) and
(\ref{eqn:d fX}),
 we have
\begin{eqnarray*}
\jetc{ \jetd X,\dM f\otimes Y}&=&\dM
f\otimes\Courant{X,Y}+\dM(a(X)f)\otimes Y,\\
\jetc{ \dM f\otimes Y,\jetd X}&=&\dM
f\otimes\Courant{Y,X}-\dM(a(X)f)\otimes Y+2\conpairing{X,Y}\jetd a^*(\dM f),\\
\jetc{\dM f\otimes X,\dM g\otimes Y}&=&a(X)(g)\dM f\otimes
Y-a(Y)(f)\dM g\otimes X.
\end{eqnarray*}

\item[c)] For any $X\in\Gamma(C)$, $\rho(\jetd X)\in\Gamma(\dev(T^*M))$
is given by
\begin{equation}\label{eqn:rho}
\rho(\jetd X)(\cdot)=\Liederivative_{a(X)}(\cdot).
\end{equation}
By  (\ref{eqn:d fX}), we get
$$
\rho(\dM f\otimes X)=a(X)\otimes \dM f,\quad \forall~f\in
C^\infty(M).
$$
For any $\xi\in\Omega^1(M)$, we have
$$\rho(\jetd X)(f\xi)=\Liederivative_{a(X)}(f\xi)=f\Liederivative_{a(X)}(\xi)+a(X)(f)\xi,$$
which implies  that $\jd\circ\rho\circ\jetd X=a(X)$, where
$\jd:\dev(T^*M)\longrightarrow TM$ is the anchor of $\dev(T^*M)$
given in (\ref{Seq:DE}). Furthermore, for any $g\in C^\infty(M)$,
the fact that $\rho(\dM f\otimes X)(g\xi)=g\rho(\dM f\otimes
X)(\xi)$ implies   that $\jd\circ\rho(\dM f\otimes X)=0$.
\end{enumerate}

We identify $C$ with $C^*$ by the bilinear form. For any $f,~g\in
C^\infty(M)$, it is straightforward to obtain the following
relations:
\begin{equation}\label{eqn:rho wx}
 \left\{\begin{array}{c}
\rhowx(\jetd\dM f)=\jetd (a^*\dM f),\\
\rhowx(\dM f\otimes \dM g)=\dM g\otimes a^*(\dM f),\\
\rhowx(\jetd( f\dM g))=\dM g\otimes a^*(\dM f)+f\jetd (a^*\dM g).
\end{array}
\right.
\end{equation}
These structures give rise to a $T^*M$-Courant algebroid.
\begin{thm}\label{thm:jet C}
For any Courant algebroid $C$, $(\jet
C,\starpair{\cdot,\cdot},\jetc{\cdot,\cdot},\rho)$ is a
$T^*M$-Courant algebroid.
\end{thm}
\pf It is straightforward to see that the pairing
$\starpair{\cdot,\cdot}$ and $\rho$ are bundle maps and
$(\Gamma(\jet C),\jetc{\cdot,\cdot}) $ is a Leibniz algebra. To show
 that the data $(\jet
C,\starpair{\cdot,\cdot},\jetc{\cdot,\cdot},\rho)$ satisfies the
properties listed  in Definition \ref{defi: $E$-C A}, it suffices to
consider elements of the form $\jetd X,~\jetd Y,~\jetd Z, ~\dM
f\otimes X,~\dM g\otimes Y,~\dM h\otimes Z$, where
$X,~Y,~Z\in\Gamma(C),~f,~g,~h\in C^\infty(M)$.

First we check Property (EC-1). Clearly, we have
$$
\rho\jetc{\jetd X,\jetd Y}=\rho\jetd\Courant{ X,
Y}=\Liederivative_{a\Courant{ X, Y}}=[\Liederivative_{a
(X)},\Liederivative_{a(Y)}]_\dev=[\rho\jetd X,\rho\jetd Y]_\dev.
$$
Furthermore, since $a\circ a^*=0$, we have
\begin{eqnarray*}
\rho\jetc{\dM f\otimes X,\jetd Y}&=&\rho(\dM
f\otimes\Courant{X,Y}-\dM(a(Y)f)\otimes X+2\conpairing{X,Y}\jetd a^*(\dM f))\\
&=&a([X,Y])\otimes\dM f-a(X)\otimes\dM(a(Y)f).
\end{eqnarray*}
On the other hand,
\begin{eqnarray*}
 {[\rho(\dM f\otimes X),\rho(\jetd
Y)]_{\dev}}(\xi)&=&[a(X)\otimes\dM f ,\Liederivative_{a(
Y)}]_{\dev}(\xi)=\langle\Liederivative_{a( Y)}\xi,a(X)\rangle\dM
f-\Liederivative_{a( Y)}(\langle
a(X),\xi\rangle\dM f)\\
&=&\langle a([X,Y]),\xi\rangle\dM f-\langle
a(X),\xi\rangle\dM(a(Y)f),
\end{eqnarray*}
which implies
$$
\rho\jetc{\dM f\otimes X,\jetd Y}=[\rho(\dM f\otimes X),\rho(\jetd
Y)]_{\dev}.
$$
Similarly, we have
\begin{eqnarray*}
\rho\jetc{\jetd X,\dM f\otimes Y}&=&[\rho(\jetd X),\rho(\dM f\otimes Y)]_\dev\\
&=&a([X,Y])\otimes\dM f+a(Y)\otimes\dM(a(X)f),\\
 \rho\jetc{\dM f\otimes X,\dM g\otimes
Y}&=&[\rho(\dM f\otimes X),\rho(\dM g\otimes
Y)]_\dev\\
&=&(a(X)g)a(Y)\otimes\dM f -(a(Y)f)a(X)\otimes\dM g.
\end{eqnarray*}
To see Property (EC-2), notice that $\jetc{\dM f\otimes X,\dM
f\otimes X}=0$ and $\starpair{\dM f\otimes X,\dM f\otimes X}=0$, so
we have
$$
\jetc{\dM f\otimes X,\dM f\otimes X}=\rhowx\jetd\starpair{\dM
f\otimes X,\dM f\otimes X}.
$$
Furthermore,
$$
\jetc{\jetd X,\jetd X}=\jetd\Courant{X,X}=\jetd\circ
a^*(\dM\conpairing{X,X})=\rhowx\jetd\circ\dM\conpairing{X,X}=\rhowx\jetd\starpair{\jetd
X,\jetd X},
$$
\begin{eqnarray*}
\jetc{\dM f\otimes X,\jetd  Y}+\jetc{\jetd  Y,\dM f\otimes
X}&=&2\conpairing{X,Y}\rhowx\jetd\dM f+2\dM f\otimes
a^*(\dM\conpairing{X,Y})\\
&=&2\rhowx\jetd(\conpairing{X,Y}\dM f)=2\rhowx\jetd\starpair{\dM
f\otimes X,Y},
\end{eqnarray*}
which implies that Property (EC-2) holds. It is straightforward to
verify Property (EC-3).
 Property (EC-4) follows from (\ref{eqn:rho wx}).   Property
(EC-5) follows from the fact that $a\circ a^*=0.$ \qed

\section{The $E$-dual pair of Lie algebroids}

Let $A$ be a vector bundle and $B$  a subbundle of $\Hom(A,E)$. For
any $\mu^k\in \Hom(\wedge^k A, E)$, denote by $\mu^{k}_\natural$ the
induced bundle map from $\wedge^{k-1}A$ to $\Hom(A, E)$ such that
\begin{equation}\label{map}
 \mu^{k}_\natural(X_1,\cdots,
X_{k-1})(X_k)=\mu^{k}(X_1,\cdots, X_{k-1},X_k).
\end{equation}
Introduce a series of vector bundles $\Hom(\wedge^k A,E)_B,~k\geq 0$
by setting $\Hom(\wedge^0 A,E)_B=E$, $\Hom(\wedge^1 A,E)_B=B$ and
\begin{eqnarray}
  \label{eqn:B k}\Hom(\wedge^k A,E)_B &\defbe& \set{\mu^{k}\in \Hom(\wedge^k A, E)~|~
\Img (\mu^{k}_\natural)\subset B}, \quad (k\geq 2).
\end{eqnarray}
If $B$ is a subbundle of $\Hom(A,E)$, then $A$ is also a bundle of
$\Hom(B,E)$. The notation $\Hom(\wedge^k B,E)_A$ is thus clear.

\begin{defi}Let $A$ and $E$ be two vector bundles over $M$. A
vector bundle $B\subset \Hom(A,E)$ is called an $E$-dual bundle of
$A$ if  the $E$-valued pairing $ \conpairing{\cdot,\cdot}_E: ~~
A\times_M B\lon E,$  ~~ $\conpairing{a,b}_E\defbe b(a)$ (where $a\in
A$, $b\in B$) is nondegenerate.
\end{defi}
Obviously, if $B$ is an $E$-dual bundle of $A$, then $A$ is an
$E$-dual bundle of $B$. We call the pair $(A,B)$ an {\bf $E$-dual
pair of vector bundles.}

Assume that $(A,[\cdot,\cdot],a)$ is a Lie algebroid and $B\subset
\Hom(A,E)$ is  an $E$-dual bundle of $A$. A representation  $\rho_A:
~A\lon \dev E$ of $A$ on $E$ is said to be {\bf{$B$-invariant}} if
$(\Gamma(\Hom(\wedge^\bullet A,E)_B),\dA)$ is a subcomplex of
$(\Omega^\bullet( A, E),\dA)$, where $\dA$ is the
  coboundary operator associated to $\rho_A$. If
$\rho_A$ is a $B$-invariant representation, we have
$\rhoAWuXing(\jet E)\subset B$. In fact, by  definition, one has
$$
\rhoAWuXing(\mu)(X) = \conpairing{\mu,\rho_A(X)}_E,\quad\forall
~\mu\in \jet E,~X\in A,
$$
and it follows that  $\rhoAWuXing: ~\jet E\lon B$  is given by
 $\rhoAWuXing([u]_m)=(\dA u)_m$, for all$ ~ u\in
\Gamma(E)$. Thus, $\rhoAWuXing(\jet E)\subset B$ is equivalent to
the condition  that $\dA(\Gamma(E))\subset \Gamma(B)$.

Furthermore, for any representation  $\rho_A: ~A\lon \dev E$, there
are  two natural  Lie derivative operations along $X\in\Gamma(A)$.
The first one is
$$\Liederivative_X:\Gamma(\Hom(\wedge^kA,E))\longrightarrow
\Gamma(\Hom(\wedge^kA,E))=\Gamma(\wedge^kA^*\otimes E)$$ defined by
\begin{eqnarray*}
\Liederivative_X(\omega\otimes u)=(\Liederivative_X\omega)\otimes
u+\omega\otimes \rho_A(X)u,\quad
\forall~\omega\in\Gamma(\wedge^kA^*),~u\in\Gamma(E).
\end{eqnarray*}
The second one is
$$\Liederivative_X:\Gamma(\Hom(\wedge^k(A^*\otimes
E),E))\longrightarrow\Gamma(\Hom(\wedge^k(A^*\otimes
E),E))=\Gamma(\wedge^k(A\otimes E^*)\otimes E)$$ defined by
$\Liederivative_X u = \rho_A(X) u$, for  $ u\in \Gamma(E)$, and
\begin{equation}
\Liederivative_X\Xi(\varpi_1\wedge\cdots\wedge\varpi_k)=\rho_A(X)(\Xi(\varpi_1\wedge\cdots\wedge\varpi_k))
-\sum_{i=1}^k\Xi(\varpi_1\wedge\cdots\wedge\Liederivative_X\varpi_i\wedge\cdots\wedge\varpi_k),
\end{equation}
for all $\Xi\in\Gamma(\Hom(\wedge^k(A^*\otimes
E),E)),~\varpi_i\in\Gamma(A^*\otimes E).$ In particular, since
$A\subset \Hom(A^*\otimes E, E)$, we have
\begin{eqnarray*}
\Liederivative_X Y=[X,Y],\quad\forall~ Y\in \Gamma(A).
\end{eqnarray*}

\begin{pro}
Let $A$ be a Lie algebroid together with a representation
$\rho:A\longrightarrow \dev E$ and $B\subset\Hom(A,E)$  a subbundle
of $\Hom(A,E)$ such that $(A,B)$ is an $E$-dual pair of vector
bundles. Then
 the following statements are equivalent:
\begin{itemize}
\item[(1)] the representation $\rho_A: ~A\lon \dev E$ is
$B$-invariant;
\item[(2)] $\dA\Gamma(E)\subset \Gamma(B)$ and
$\dA\Gamma(B)\subset \Gamma(\Hom(\wedge^2 A,E)_B)$;
\item[(3)] $\Gamma(\Hom(\wedge^k A,E)_B)$ is
 invariant under the operation  $\Liederivative_{X}$ for any $X\in\Gamma(A)$;
\item[(4)] $\Gamma(\Hom(\wedge^k B,E)_A)$  is
 invariant under the operation $\Liederivative_{X}$ for any $X\in\Gamma(A)$.
\end{itemize}
\end{pro}
\pf The implication   $(1)\Longrightarrow (2)$ is obvious. We adopt
an inductive approach to see the implication   $(2)\Longrightarrow
(1)$.  For any $n\geq 1$, $\Liederivative_{X}:\Gamma(\Hom(\wedge^n
A,E)_B)\longrightarrow\Gamma(\Hom(\wedge^n A,E)_B)$ is  well defined
and we have $i_{Y} \Liederivative_X -\Liederivative_X
i_{Y}=i_{[Y,X]}$. Assume that $\dA \Gamma( \Hom(\wedge^{n-1}
A,E)_B)\subset \Gamma(\Hom(\wedge^n A,E)_B)$ and $\dA \Gamma
(\Hom(\wedge^n A,E)_B)\subset \Gamma(\Hom(\wedge^{n+1} A,E)_B)$ hold
 for all $\mu^{n+1}\in
\Gamma( \Hom(\wedge^{n+1} A,E)_B)$. To prove that $\dA \mu^{n+1}\in
\Gamma(\Hom(\wedge^{n+2} A,E)_B)$, it suffices to show that $i_X\dA
\mu^{n+1}\in \Gamma(\Hom(\wedge^{n+1} A,E)_B)$, for all $a\in
\Gamma(A)$. Again, it suffices to show that $i_{Y}i_X\dA
\mu^{n+1}\in \Gamma(\Hom(\wedge^n A,E)_B)$ holds for all $Y\in
\Gamma(A)$. In fact,
\begin{eqnarray*}
i_{Y}i_X\dA \mu^{n+1}&=&i_{Y}(\Liederivative_X \mu^{n+1}-\dA i_X \mu^{n+1})\\
&=&(i_{Y} \Liederivative_X  -\Liederivative_X i_{Y})\mu^{n+1} +
\Liederivative_Xi_{Y}\mu^{n+1}-i_{Y}\dA i_X \mu^{n+1}\\
&=& i_{[Y,X]}\mu^{n+1}+ \Liederivative_Xi_{Y}\mu^{n+1}-i_{Y}\dA i_X
\mu^{n+1}\in \Gamma(\Hom(\wedge^n A,E)_B).
\end{eqnarray*}
So we conclude that $\Gamma(\Hom(\wedge^\bullet A,E)_B)$ is a
subcomplex of $\Omega^\bullet( A, E)$. This completes the proof of
the equivalence of $(1)$ and $(2)$. The equivalence of $(1)$ and
$(3)$ is obvious.

Next we prove the  equivalence of $(2)$ and $(4)$. For any
$X^k\in\Gamma(\Hom(\wedge^k B,E)_A) $ and   $\xi_i\in B$, we have
\begin{eqnarray*}
&&\conpairing{i_{\xi_1\wedge\cdots\wedge \xi_{k-1}}\Liederivative_X
X^k,\xi_k}_E\\&=&(\Liederivative_X X^k)(\xi_1\wedge \xi_2\wedge
\cdots
\wedge \xi_k)\\
& =& \rho_A(X) (X^k (\xi_1\wedge \xi_2\wedge \cdots \wedge
\xi_k))-\sum_{i=1}^k X^k(\xi_1\wedge\cdots \wedge\Liederivative_X
\xi_i
\wedge \cdots \wedge \xi_k)\\
 &=&\rho_A(X)\conpairing{i_{\xi_1\wedge\cdots\wedge \xi_{k-1}}X^k,\xi_k}_E-
 \sum_{j=1}^{k-1}\conpairing{i_{\xi_1\wedge\cdots\wedge\Liederivative_X\xi_j\wedge\cdots\wedge \xi_{k-1}}X^k,\xi_k}_E
 -\conpairing{i_{\xi_1\wedge\cdots\wedge
 \xi_{k-1}}X^k,\Liederivative_X\xi_k}_E\\
 &=&\conpairing{[X,i_{\xi_1\wedge\cdots\wedge \xi_{k-1}}X^k]-\sum_{j=1}^{k-1}i_{\xi_1\wedge\cdots\wedge\Liederivative_X\xi_j\wedge\cdots\wedge
 \xi_{k-1}}X^k,\xi_k}_E.
\end{eqnarray*}
Since the $E$-valued pairing $\langle\cdot,\cdot\rangle_E$ is
nondegenerate, we have
$$
i_{\xi_1\wedge\cdots\wedge \xi_{k-1}}\Liederivative_X
X^k=[X,i_{\xi_1\wedge\cdots\wedge
\xi_{k-1}}X^k]-\sum_{j=1}^{k-1}i_{\xi_1\wedge\cdots\wedge\Liederivative_X\xi_j\wedge\cdots\wedge
 \xi_{k-1}}X^k,
$$
which implies the equivalence of $(2)$ and $(4)$.  \qed

\begin{defi}
An \textbf{$E$-dual pair of Lie algebroids}
$((A,\rho_A);(B,\rho_B))$ consists of two Lie algebroids $A$ and $B$
which are mutually $E$-dual vector bundles, a $B$-invariant
representation $\rho_A: ~A\lon \dev E$ and an $A$-invariant
representation $\rho_B: ~B\lon \dev E$.
\end{defi}

Obviously,  $\jet E$ and $\dev E$ are mutually  $E$-dual bundles. In
the following, we show some properties of the  bundle $
\Hom(\wedge^k \dev E,E)_{\jet E}$.

\begin{pro}\label{Pro:muk=phi gamma}
If $k\geq 2$, for any $\mu^k\in \Gamma(\Hom(\wedge^k \dev E,E)_{\jet
E})$, there is a unique bundle map $\lambda_{\mu^k}\in
\Gamma(\Hom(\wedge^{k-1}TM, E))$ such that
$$ \mu^k( \frkd_1\wedge\cdots\wedge\frkd_{k-1}\wedge \Phi)=
\Phi\circ
{\lambda_{\mu^k}}(\jd(\frkd_1)\wedge\cdots\wedge\jd(\frkd_{k-1})),\
\quad\forall ~~ \Phi\in \Gamma(\gl(E)), \ \frkd_i\in
\Gamma(\dev{E}).$$
\end{pro}
\pf By definition, we have $
\mu^{k}(\frkd_1\wedge\cdots\wedge\frkd_k)=\conpairing{\mu^{k}_\natural(\frkd_1\wedge\cdots\wedge\frkd_{k-1}),\frkd_k}_E,
$ which implies
$$
\p\circ
\mu^{k}_\natural(\frkd_1\wedge\cdots\wedge\frkd_{k-1})=\mu^k(\frkd_1\wedge\cdots\wedge\frkd_{k-1}\wedge
\Id).
$$
We claim that $i_{\Phi}(\p\circ\mu^{k}_\natural)=0$, for all $~
\Phi\in \Gamma(\gl(E))$. In fact, we have
\begin{eqnarray*}
&&\p\circ
\mu^{k}_\natural(\frkd_1\wedge\cdots\wedge\frkd_{k-2}\wedge \Phi)\\
&=&\mu^k(\frkd_1\wedge\cdots\wedge\frkd_{k-2}\wedge \Phi\wedge
\Id)=-\mu^k(\frkd_1\wedge\cdots\wedge\frkd_{k-2}\wedge
\Id\wedge \Phi  )\\
&=&-\conpairing{\mu^{k}_\natural(\frkd_1\wedge\cdots\wedge\frkd_{k-2}\wedge
\Id),\Phi}_E =-\Phi\circ \p\circ
\mu^{k}_\natural(\frkd_1\wedge\cdots\wedge\frkd_{k-2}\wedge
\Id)\\&=&-\Phi\circ \mu^k(\frkd_1\wedge\cdots\wedge\frkd_{k-2}\wedge
\Id\wedge \Id)\\&=&0.
\end{eqnarray*}
Therefore, the bundle map $\p\circ \mu^{k}_\natural: \wedge^{k-1}
\Gamma(\dev E)\lon \Gamma(E)$ factors through $\jd$, i.e. there is a
unique bundle map ${\lambda_{\mu^k}}: \Gamma(\wedge^{k-1}TM)\lon
\Gamma(E)$ such that
$$\p\circ
\mu^{k}_\natural(\frkd_1\wedge\cdots\wedge\frkd_{k-1})=
{\lambda_{\mu^k}}(\jd(\frkd_1)\wedge\cdots\wedge\jd(\frkd_{k-1})).\qed$$
   \vspace{2mm}

Therefore, for $k\geq2$, the vector bundle $\Hom(\wedge^k \dev
E,E)_{\jet E}$ can be defined directly by
\begin{eqnarray*}
&& \Hom(\wedge^k \dev E,E)_{\jet E} \\&=&  \{\mu^k \in
\Hom(\wedge^k\dev{E},E)~~|~~\exists~!~ {\lambda_{\mu^k}}\in
\Hom(\wedge^{k-1}TM, E),\mbox{ s.t., }~\forall ~~ \Phi\in \gl(E),
\\&& \frkd_i\in \dev{E},
\qquad\mu^k(\frkd_1\wedge\cdots\wedge\frkd_{k-1}\wedge \Phi)=
\Phi\circ
{\lambda_{\mu^k}}(\jd(\frkd_1)\wedge\cdots\wedge\jd(\frkd_{k-1}))\}.
\end{eqnarray*}
We will write ${\lambda_{\mu^k}}={\p}^k(\mu^k)$ for $\mu^k$ given
above. For any $\xi\in\Hom(\wedge^k TM,E)$, $k\geq1$, we define
${\e}^k(\xi)\in \Hom(\wedge^k \dev E,E)_{\jet E}$ by
$$
{\e}^k(\xi)(\frkd_1\wedge\cdots\wedge\frkd_k)\defbe
\xi{(\jd(\frkd_1)\wedge\cdots\wedge\jd(\frkd_k))},\quad \forall ~~
\frkd_i\in \dev{E}.
$$
In addition, we regard $\Hom(\wedge^{-1}TM,E)=0$,
$\Hom(\wedge^0TM,E)=E$ and $\Hom(\wedge^0 \dev E,E)_{\jet E}=E$. Let
$\p^0=0$ and $\e^0=\Id_E$.
\begin{pro}
For any $k\geq0$, the following sequence is exact:
\begin{equation}\label{temp121adfa2}
0 \longrightarrow \Hom(\wedge^k TM,
E)\stackrel{{\e}^k}{\longrightarrow}
                  {\Hom(\wedge^k \dev E,E)_{\jet E}} \stackrel{{\p}^k}{\longrightarrow} \Hom(\wedge^{k-1}TM, E) \longrightarrow
                  0.
\end{equation}
\end{pro}
\pf If $k=0,1$, the result is clear. For $k\geq2$,  it is obvious
that ${\e}^k$ is an injection and ${\p}^k\circ {\e}^k=0$. Now we
prove that ${\p}^k$ is  surjective. For every ${\lambda}\in
\Gamma(\Hom(\wedge^{k-1}TM, E))$, i.e. a bundle map from
$\wedge^{k-1}TM$ to $E$, we define $\widetilde{{\lambda}}\in
\Gamma(\Hom(\wedge^k \dev E,E)_{\jet E})$ by
\begin{eqnarray*}
\widetilde{{\lambda}} (\frkd_1\wedge\cdots\wedge\frkd_{k }) &\defbe
& \sum_{i=1}^{k } (-1)^{i+1} \frkd_i\circ
{\lambda}(\jd(\frkd_1)\wedge\cdots\widehat{\jd(\frkd_i)}\cdots\wedge\jd(\frkd_{k
}) )\\
&&
+\sum_{i<j}(-1)^{i+j}\lambda([\jd(\frkd_i),\jd(\frkd_j)]\wedge\jd(\frkd_1)\wedge\cdots\widehat{\jd(\frkd_i)}\cdots\widehat{\jd(\frkd_j)}\cdots\wedge\jd(\frkd_{k})),
\end{eqnarray*}
for any $\frkd_i\in\Gamma(\dev E),~i=1,\cdots,k. $ It is
straightforward to see that $\widetilde{\lambda}$ is a bundle map
from $\wedge^k \dev E$ to $E$. Since for any
$\Phi\in\Gamma(\gl(E))$,
$$
\widetilde{{\lambda}}( \frkd_1\wedge\cdots\wedge\frkd_{k-1}\wedge
\Phi)= (-1)^{k+1}\Phi\circ
{\lambda}(\jd(\frkd_1)\wedge\cdots\wedge\jd(\frkd_{k-1})),
$$
we have $\widetilde{{\lambda}}\in\Gamma(\Hom(\wedge^k \dev
E,E)_{\jet E})$ and $\p^k
\widetilde{{\lambda}}=(-1)^{k+1}{\lambda}$.

Finally, if $\mu^k\in \Gamma(\Hom(\wedge^k \dev E,E)_{\jet E})$
satisfies $ {\p}^k(\mu^k)=0$, then we have
$$
\mu^k(\frkd_1\wedge\cdots\wedge \frkd_{k-1}\wedge \Phi)= \Phi\circ
{\p}^k(\mu^k) (\jd(\frkd_1)\wedge \cdots\wedge \jd(\frkd_{k-1}))=0,
$$
which implies that the map $\mu^k $ factors through $\jd$, i.e.
there is a unique $\xi\in \Hom(\wedge^{k}TM, E)$ such that
$$
\mu^k(\frkd_1\wedge\cdots\wedge \frkd_{k }
 )= \xi
 (\jd(\frkd_1)\wedge\cdots\wedge\jd(\frkd_k)) .
$$
Therefore,  sequence (\ref{temp121adfa2}) is exact. \qed
\vspace{3mm}

 In the sequel, we will omit the embedding ${\e}^k$ and
directly regard $\Hom(\wedge^k TM, E)$  as a subbundle  of $
\Hom(\wedge^k \dev E,E)_{\jet E}$. Hence we have
\begin{equation}\label{eqn:property p i}
\p(i_\frkd\mu^k)=i_{\jd(\frkd)}\p^k\mu^k,\quad\forall~\frkd\in\dev
E,~\mu^k\in\Hom(\wedge^k \dev E,E)_{\jet E}.
\end{equation}

Recall that $\dev{E}$ is a Lie algebroid and there is a natural
representation $\Id_{\dev E}$ on $E$. For $\mu^k\in
\Gamma(\Hom(\wedge^k \dev E,E)_{\jet E})$ and $\frkd_i\in\Gamma(\dev
E)$, the coboundary operator $\jetd: \Omega^\bullet(\dev E, E)\lon
\Omega^{\bullet+1}(\dev E, E)$
 is given by
\begin{eqnarray}\nonumber
\jetd \mu^k (\frkd_1\wedge\cdots\wedge\frkd_{k+1}) &\defbe &
\sum_{i=1}^{k+1} (-1)^{i+1} \frkd_i\circ
\mu^k(\frkd_1\wedge\cdots\widehat{\frkd_i}\cdots\wedge\frkd_{k+1}
)\\\label{Eqt:jetd} &&
+\sum_{i<j}(-1)^{i+j}\mu^k([\frkd_i,\frkd_j]_{\dev}\wedge\frkd_1\wedge\cdots\widehat{\frkd_i}\cdots\widehat{\frkd_j}\cdots\wedge\frkd_{k+1}).
\end{eqnarray}

\begin{lem}\label{Lem:property p i} The representation $\Id_{\dev E}$ of the gauge Lie algebroid $\dev E$ on $E$ is $\jet E$-invariant,
i.e. $(\Gamma(\Hom(\wedge^\bullet\dev E, E)_{\jet E}),\jetd)$ is a
subcomplex of $(\Gamma(\Hom(\wedge^\bullet\dev E, E)),\jetd)$. More
precisely, for any $\frkd\in \Gamma (\dev E)$
 and $\mu^k\in\Gamma(\Hom(\wedge^k\dev E,E)_{\jet E}) $, $k\geq0$, we have
\begin{eqnarray}
\label{Eqn:p L d mu k}
{\p}^k(\Liederivative_{\frkd}\mu^k)&=&\Liederivative_{\frkd}
({\p}^k(\mu^k)),\\
\label{lem:property p 1}\p^k\jetd \p^k\mu^k&=&(-1)^{k+1}\p^k\mu^k,\\
\label{lem:proerty p}\p^{k+1}\jetd \mu^k &=&\jetd
\p^k\mu^k+(-1)^k\mu^k.
\end{eqnarray}
\end{lem}
\pf  Assume that  $\p^k\mu^k={\lambda}\in
\Gamma(\Hom(\wedge^{k-1}TM,E))$, i.e. for all $~~\Phi\in
\Gamma(\gl(E)),~\frkd_i\in\Gamma(\dev E)$,
$$
\mu^k( \frkd_1\wedge\cdots\wedge\frkd_{k-1}\wedge \Phi)= \Phi\circ
{\lambda}(\jd(\frkd_1)\wedge\cdots\wedge\jd(\frkd_{k-1})).
$$
Then, we have
\begin{eqnarray*}
&&\Liederivative_{\frkd}\mu^k(\frkd_1\wedge\cdots\wedge\frkd_{k-1}\wedge\Phi)\\&=&\frkd\circ\mu^k(\frkd_1\wedge\cdots\wedge\frkd_{k-1}\wedge\Phi)
-\mu^k(\frkd_1\wedge\cdots\wedge\frkd_{k-1}\wedge
[\frkd,\Phi])\\&&-\sum_{i=1}^{k-1}\mu^k(\frkd_1\wedge\cdots\wedge[\frkd,
\frkd_i]_{\dev}\wedge\cdots\wedge\frkd_{k-1}\wedge\Phi)\\
&=&\frkd\circ\Phi\circ\lambda(\jd\frkd_1\wedge\cdots\wedge\jd\frkd_{k-1})-[\frkd,\Phi]\circ\lambda(\jd\frkd_1\wedge\cdots\wedge\jd\frkd_{k-1})\\
&&-\sum_{i=1}^{k-1}\Phi\circ\lambda(\jd\frkd_1\wedge\cdots\wedge[\jd\frkd,\jd\frkd_i]\wedge\cdots\wedge\jd\frkd_{k-1})\\
&=&\Phi(\frkd\circ\lambda(\jd\frkd_1\wedge\cdots\wedge\jd\frkd_{k-1})-\sum_{i=1}^{k-1}\lambda(\jd\frkd_1\wedge\cdots\wedge[\jd\frkd,\jd\frkd_i]\wedge\cdots\wedge\jd\frkd_{k-1}))\\
&=&\Phi(\Liederivative_{\frkd}\lambda(\frkd_1\wedge\cdots\wedge\frkd_{k-1})),
\end{eqnarray*}
which implies   the equality
${\p}^k(\Liederivative_{\frkd}\mu^k)=\Liederivative_{\frkd}
({\p}^k(\mu^k))$. The other conclusions can be obtained similarly.
\qed

\begin{rmk}
For any ${\lambda}\in
\Gamma(\Hom(\wedge^{k-1}TM,E))\subset\Gamma(\Hom(\wedge^{k-1}\dev
E,E)_{\jet E})$, $\Liederivative_{\frkd}\lambda$ can be considered
as the Lie derivative of $ \Omega^{k-1}(M)\otimes \Gamma(E)$ in the
obvious sense:
$$
\Liederivative_{\frkd}\defbe \Liederivative_{\jd( \frkd)}\otimes
\Id_{\Gamma(E)}+\Id_{\Omega^{k-1}(M)}\otimes \frkd.
$$
\end{rmk}

By Lemma \ref{Lem:property p i}, the representation $\Id_{\dev E}$
of $\dev E$ on $E$ is $\jet E$-invariant. Since $\jet E$ is a Lie
algebroid with all  structures zero, we have
\begin{cor}The pair
$((\dev E,\Id_{\dev E});(\jet E,0))$ is an $E$-dual pair.
\end{cor}
\begin{thm}\label{Thm:longExact}
For the cochain complex $C(E)=(\Gamma(\Hom(\wedge^\bullet\dev
E,E)_{\jet E}),\jetd)$, we have $\mathrm{H}^k(C(E))=0$, for all
$k=0,1,2,\cdots$. In other words, there is a long exact sequence:
\begin{eqnarray}
\nonumber0\lon\Gamma(E)\stackrel{\jetd}{\longrightarrow}\Gamma(\jet{E})
\stackrel{\jetd}{\longrightarrow} \Gamma(\Hom(\wedge^2\dev
E,E)_{\jet
E})\stackrel{\jetd}{\longrightarrow}\cdots\\\label{seq:long
exact}\cdots\stackrel{\jetd}{\longrightarrow}
\Gamma(\Hom(\wedge^n\dev E,E)_{\jet E})\lon 0,
\end{eqnarray}
where $n=\dimension M+1$.
\end{thm}
\pf If $k=\dimension M+2$, we have
$\Hom(\wedge^{k+2}TM,E)=\Hom(\wedge^{k+1}TM,E)=0$. By the exact
sequence (\ref{temp121adfa2}),  $\Hom(\wedge^{k+2}\dev E,E)_{\jet
E}=0$. By (\ref{lem:proerty p}) and if $\jetd\mu^k=0$, we have
$\mu^k=(-1)^{k+1}\jetd(\p^k\mu^k)$, and hence
$\mathrm{H}^k(C(E))=0$. \qed


\begin{ex}{\rm
If $E$ reduces to a vector space $V$, $\jet V=V$, $\dev V=\gl(V)$.
By (\ref{temp121adfa2}), we have $\Hom(\wedge^2\gl(V),V)_V=0$. In
fact, for any $\phi\in\Hom(\wedge^2\gl(V),V)_V$, $~A,~B\in\gl(V)$,
we have
$$
\phi(A\wedge B)=B\phi(A\wedge  \Id_V)=-BA\phi(\Id_V\wedge \Id_V)=0.
$$
Therefore, $\phi=0$, which implies that
$\Hom(\wedge^2\gl(V),V)_V=0$. On the other hand, for any $v\in V$,
$\jetd v=v$ and for any $u\in \jet V$, $\jetd u=0$. Thus, the first
cohomology is trivial.}
\end{ex}
\begin{ex}\rm{
Let $E=M\times \mathbb R$, then we have $\dev E=TM\oplus \mathbb R$
and $\jet E=T^*M\oplus \mathbb R$. Denote the basis of
$C^\infty(M)\subset \Gamma(\dev E)$ by $\Id$ and the dual basis by
$\Id^*$, i.e. $\Id^*(f\Id)=f$, for all $f\in C^\infty(M)$. Since
$\rho:\dev E\oplus\jet E\longrightarrow \dev E$ is the projection,
we can write $\rho=pr_{TM}+\Id^*$, where $\Id^*$ is considered as a
section of $\omni^*$ satisfying
$\Id^*(f)=f,~\Id^*(X)=\Id^*(\mu)=0,~$ for all $f\in
C^\infty(M),~X\in \frkX(M) $ and $\mu\in\Gamma(\jet E)$. Since
$\jetd f=\dM f+f\Id^*$,  we have
$$
\jetd f=0\Longleftrightarrow f=0,\quad \forall ~f\in C^\infty(M),
$$
 which implies that $\mathrm{H}^0(C(E))=0$.

 For any $\lambda\in \Omega^1(M)$, by
(\ref{lem:proerty p}), one gets $\p^2(\jetd \lambda)=-\lambda$.
Furthermore, we have $\jetd (f\Id^*)=\jetd (\jetd(f)-\dM
f)=-\jetd(\dM f)$, which implies that $\p^2(\jetd (f\Id^*))=\dM f$.
Therefore, for any $\mu\in \Gamma(\jet E)$, we have
$$\jetd
\mu=0\Longleftrightarrow \mu=\dM f+f\Id^*\Longleftrightarrow
\mu=\jetd f,
$$
for some $f\in C^\infty(M)$, which implies that
$\mathrm{H}^1(C(E))=0$. For similar reasons, one has
$\mathrm{H}^k(C(E))=0$.}
\end{ex}
\section{The automorphism groups of  omni-Lie algebroids  }
In this section, we study the automorphism groups and the twists of
omni-Lie algebroids. For $i=1,2$, it is subtle to define morphisms
between two $E_i$-Courant algebroids  with different base manifolds,
which remains a topic in the future.
 As for general Lie algebroid
morphisms and   Courant algebroid morphisms, please refer to
\cite{CL2,PJHK2,Courant morphism}.   Here we only consider the
automorphisms of $E$-Courant algebroids.

 Given an automorphism
$\Phi:E\longrightarrow
 E$ over  the diffeomorphism  $\phi:~M\lon M$ of the base manifold,
 there  induces a unique automorphism $\Ad_\Phi$ of
 $\dev E$ such that:
 $$\Ad_\Phi(\frkd)(u)=\Phi\circ \frkd\circ\Phi^{-1}(u),\quad
 \forall~\frkd\in \Gamma(\dev E),~u\in\Gamma(E).
 $$
\begin{defi}
The  automorphism group $\Aut(\huaK)$ of an $E$-Courant algebroid
$\huaK$ is the group of bundle automorphisms
$F:\huaK\longrightarrow\huaK$ covering bundle automorphisms
$\Phi:E\longrightarrow E$ such that
\begin{itemize}
\item[\rm(1) ]$F$ is orthogonal, i.e. $\ppairingE{F(X),F(Y)}=\Phi\ppairingE{X,Y}$;
\item[\rm(2) ]$F$ is bracket-preserving, i.e.
$F\CDorfman{X,Y}=\CDorfman{F(X),F(Y)}$;
\item[\rm(3) ]$F$ is
compatible with the anchor, i.e. $\rho\circ F=\Ad_\Phi\circ\rho$.
\end{itemize}
\end{defi}

We usually denote such an automorphism  by a pair $(F,\Phi)$. The
set of all automorphisms $( F, \Id_E)$ is a normal subgroup of
$\Aut(\huaK)$,  similar to the $B$-field introduced in
\cite{GualtieriGeneralizedComplex}.

 Now let us  study the automorphism
group of the omni-Lie algebroid $\omni=\dev E\oplus\jet E$ defined
in Definition \ref{def:omni algebroid}. For any automorphism
$\Phi:E\longrightarrow
 E$, there is an induced map $\widetilde{\Phi}:\jet E\longrightarrow\jet
 E$ defined by
 $$
\widetilde{\Phi}(\mu)=[\Phi(u)]_{\phi(m)},\quad
\forall~\mu=[u]_m\in(\jet E)_m,~u\in\Gamma(E).
$$
It is clear that the pair $(\Ad_\Phi+\widetilde{\Phi},\Phi)$ is an
 automorphism of the omni-Lie algebroid $\omni$ and it is
totally determined by  $\Phi$.

There is another symmetry of the omni-Lie algebroid $\omni$, which
we call the $B$-field transformation. Let us elaborate on this idea.
For any $b\in \Gamma(\Hom(\dev E,\jet E))$, there is a
transformation $e^b: \omni\lon \omni$ defined by
\begin{equation}\label{e b}
e^b\left(\begin{array}{c}\frkd\\\mu\end{array}\right)
=\left(\begin{array}{cc}\Id&0\\b&\Id\end{array}\right)\left(\begin{array}{c}\frkd\\\mu\end{array}\right)=\left(\begin{array}{c}\frkd\\\mu+i_\frkd
b\end{array}\right).
\end{equation}

\begin{lem}
For $b\in \Gamma(\Hom(\wedge^2\dev E,E)_{\jet E})$, the map $e^b$ is
an automorphism of the omni-Lie algebroid $\omni$ if and only if $b$
is closed, i.e. $\jetd b=0$.
\end{lem}
\pf Let $\frkd,~\frkr\in\Gamma(\dev E),~\mu,~\nu\in\Gamma(\jet{E})$.
First, $b$ is skew-symmetric implies that $e^b$ preserves the
standard pairing given in (\ref{standard pair}). We also have
\begin{eqnarray}
\nonumber\Dorfman{e^b(\frkd+\mu),e^b(\frkr+\nu)}&=&\Dorfman{\frkd+\mu,\frkr+\nu}+\Dorfman{\frkd,i_\frkr
 b}+\Dorfman{i_\frkd b,\frkr }\\
\nonumber &=&\Dorfman{\frkd+\mu,\frkr+\nu}+\Liederivative_\frkd
i_\frkr
 b-i_\frkr\jetd i_\frkd b\\
 \label{eqn: e B}&=&e^b(\Dorfman{\frkd+\mu,\frkr+\nu})+i_\frkr i_\frkd\jetd b.
\end{eqnarray}
So $e^b$ is an automorphism of the omni-Lie algebroid $\omni$ if and
only if $i_\frkr i_\frkd\jetd b=0$ for all
$\frkd,~\frkr\in\Gamma(\dev E)$, which happens if and only if
$\jetd b=0$. \qed \vspace{3mm}

The transformation $e^b$ defined by (\ref{e b}) will be called a
{\bf{$B$-field transformation},} for any $b\in
\Gamma(\Hom(\wedge^2\dev E,E)_{\jet E})$ with $\jetd b=0$.
\begin{cor}
The abelian group of  $B$-field transformations is isomorphic to
$\Gamma(\Hom(TM, E))$.
\end{cor}
\pf For any $b\in \Gamma(\Hom(\wedge^2\dev E,E)_{\jet E})$ such that
$\jetd b=0$, assume that $\mu=\p^2b\in\Gamma(\Hom(TM,E))$, where
$\p^k$ is given in Sequence (\ref{temp121adfa2}). Then, by
(\ref{lem:proerty p}), $b=-\jetd \mu$. As a vector space,
$\Gamma(\jet E)\cong\Gamma(\Hom(TM, E))\oplus\jetd\Gamma(E)$. Since
$\jetd^2=0$, it follows that
$$\jetd(\Gamma(\jet E))\cong\jetd(\Gamma(\Hom(TM,
E)))\cong\Gamma(\Hom(TM, E)). \qed$$

In fact, any automorphism of the omni-Lie algebroid $\omni$ is a
composition of an automorphism $\Phi$ of the vector bundle $E$ and a
 $B$-field transformation.

\begin{thm}\label{thm:aut omni}
Let $(F,\Phi)$ be an  automorphism of the omni-Lie algebroid
$\omni$, where $\Phi$ is an automorphism of $E$ and
$F:\omni\longrightarrow\omni$ is an  automorphism of $\omni$. Then
$F$ can be decomposed as a composition of an automorphism $\Phi$ of
the vector bundle $E$ and a $B$-field transformation $e^b$.
\end{thm}
\pf Since $\Phi$ is an automorphism of $E$,
$(\Ad_\Phi+\widetilde{\Phi},\Phi)$ is an  automorphism of the
omni-Lie algebroid $\omni$. Setting $G=\Phi^{-1}\circ F$, the pair
$(G,\Id_E)$ is again an  automorphism of the omni-Lie algebroid
$\omni$.
Since $G$ and $\rho$ are compatible, we can write
$$G(\frkd+\mu)=\frkd+b(\frkd)+\sigma(\mu),\quad\forall ~\frkd+\mu\in\omni,$$
where $b: ~\dev E\lon \jet E$ and $\sigma: \jet E\lon \jet E$ are
two bundle maps. Then, by
$$
\ppairingE{G(\frkd+\mu),G(\frkr+\nu)}=\ppairingE{ \frkd+\mu ,
\frkr+\nu },\quad\forall ~\frkd+\mu,~\frkr+\nu\in\omni,
$$
we know that $\sigma=\Id_{\jet E}$ and $b$ is skew-symmetric:
$$
\conpairing{b(\frkd),\frkr}_E=-\conpairing{b(\frkr),\frkd}_E\,.\quad\forall~
\frkd,\frkr\in \dev E.
$$
Using the equation
$$
\Dorfman{G(\frkd+\mu),G(\frkr+\nu)}=G\Dorfman{ \frkd+\mu , \frkr+\nu
},\quad\forall ~\frkd+\mu,~\frkr+\nu\in\Gamma(\omni),
$$
we see that $b$ is closed with respect to the Lie algebroid
cohomology of $\dev E$. Thus  $F=\Phi\circ e^b$,  as required. \qed
\begin{cor}
The automorphism group $\Aut(\omni)$ of the omni-Lie algebroid
$\huaE$ is the semidirect product of $\Aut(E)$ and $\Gamma(\Hom(TM,
E))$, i.e.
\begin{equation}\label{eqn:auto}
\Aut(\omni)\cong\Aut(E)\ltimes\Gamma(\Hom(TM, E)),
\end{equation}
where the action of $\Aut(E)$ on $\Gamma(\Hom(TM, E))$,  denoted by
$\cdot$, is given by
$$
\Phi\cdot \eta=\Phi\circ\eta\circ\phi_*^{-1},$$ where $\phi_*$ is
the tangent   of the  map $\phi$ induced  by $\Phi$ on the base
manifold $M$.
\end{cor}

Differentiating a 1-parameter family of automorphisms
$F_t=\Phi_t\circ e^{tb},~F_0=\Id,~b=-\jetd\mu$, we see that the Lie
algebra $\Der(\omni)$ of infinitesimal symmetries of the omni-Lie
algebroid $\omni$ consists of  pairs $(\frkd,\mu)\in\Gamma(\dev
E)\oplus\Gamma(\Hom(TM, E))$. The pair $(\frkd,\mu)$ acts on
$\Gamma(\omni)$ via
\begin{equation}\label{eqn:der}
(\frkd+\mu)\cdot(\frkr+\nu)=[\frkd,\frkr]_\dev+\Liederivative_\frkd\nu-i_\frkr\jetd
\mu,\quad\forall~\frkr+\nu\in\Gamma(\omni).
\end{equation}

From Theorem \ref{thm:aut omni}, we  conclude:
\begin{pro}
The Lie algebra $\Der(\omni)$ of infinitesimal symmetries of an
omni-Lie algebroid $\omni=\dev{E}\oplus \jet{E}$ is isomorphic to
the semidirect sum of $\Gamma(\dev E)$ and $\Gamma(\Hom(TM, E))$,
i.e.
$$
\Der(\omni)\cong\Gamma(\dev E)\ltimes\Gamma(\Hom(TM, E)).
$$
Moreover, all of these  derivations are defined by the standard
bracket (\ref{standard bracket}) of the omni-Lie algebroid $\omni$
from the left hand side, and there is an exact sequence of Leibniz
algebra morphism:
\begin{equation}\label{seq:der}
\xymatrix@C=0.5cm{0 \ar[r] & \Gamma( E)  \ar[rr]^{\jetd} &&
                \Gamma(\omni)  \ar[rr]^{\mathrm{ad}} && \Der( \omni) \ar[r]  & 0.
                }
\end{equation}
\end{pro}
\pf By (\ref{eqn:der}), for any derivation $\frkd+\mu$, we have
$$(\frkd+\mu)\cdot(\frkr+\nu)=\Dorfman{\frkd+\mu,\frkr+\nu},$$
which implies that  derivations of the omni-Lie algebroid $\omni$
are defined by the standard bracket (\ref{standard bracket}). Also
$\jetd^2=0$ implies that   $\jetd(\Gamma(E))$ is the left center of
the standard bracket (\ref{standard bracket}), i.e. the kernel of
the map $\mathrm{ad}$. Thus, Sequence (\ref{seq:der}) is exact. \qed
\vspace{3mm}


 Similar to the fact that  an exact Courant algebroid
can be twisted by a closed 3-form,  we consider the deformation of
the omni-Lie algebroid $\omni$. Given a linear map
$\Theta:\Gamma(\dev E\otimes\dev{E})\longrightarrow\Gamma(\jet E)$,
we define a new bracket $\Dorfman{\cdot,\cdot}_\Theta$ on
$\Gamma(\dev E\oplus\jet E)$ by
$$
\Dorfman{ \frkd+\mu , \frkr+\nu}_\Theta=\Dorfman{ \frkd+\mu ,
\frkr+\nu}+\Theta(\frkd\otimes\frkr).
$$
To meet Property (EC-2), $\Theta $ must be skew-symmetric. To
satisfy Property (EC-3), we need
\begin{eqnarray*}
0&=&\ppairingE{\Dorfman{\frkd,\frkr}_\Theta ,\frkt}
+\ppairingE{\frkr,\Dorfman{\frkd,\frkt}_\Theta }
=\frac{1}{2}(\pairing{\Theta(\frkd\wedge
\frkr),\frkt}_E+\pairing{\Theta(\frkd\wedge\frkt),\frkr}_E)\\
&=&\frac{1}{2}(\pairing{\Theta(\frkd\wedge
\frkr),\frkt}_E-\pairing{\Theta(\frkt\wedge\frkd),\frkr}_E).
\end{eqnarray*}
Thus,
\begin{equation}\label{Eqt:UpsilonSkewsymmetric}
\conpairing{\Theta({\frkd}\wedge {\frkr}),
\frkt}_E=\conpairing{\Theta({\frkr}\wedge {\frkt}), \frkd}_E
=\conpairing{\Theta({\frkt}\wedge {\frkd}), \frkr}_E, \quad\forall
~~ \frkd,\frkr,\frkt\in\dev E,
\end{equation}
which implies  that $\Theta\in \Gamma(\Hom(\wedge^3\dev E,E)_{\jet
E})$. In this way, it is standard to prove that
$\Dorfman{\cdot,\cdot}_\Theta$ defines a new $E$-Courant algebroid
structure on $\dev E\oplus\jet E$ (using the standard pairing
(\ref{standard pair}) and the same anchor  of the omni-Lie
algebroid) if and only if $\jetd \Theta=0$. We call this $E$-Courant
algebroid the {\bf{$\Theta$-twisted omni-Lie algebroid}}.

\begin{thm}\label{thm:isomorphic to omni}
Any twisted omni-Lie algebroid  is isomorphic to the standard
omni-Lie algebroid $\omni$ in Definition \ref{def:omni algebroid}.
\end{thm}
\pf By Theorem \ref{Thm:longExact}, there is some
$b\in\Gamma(\Hom(\wedge^2\dev E,E)_{\jet E})$ such that
$\Theta=\jetd b$. By (\ref{eqn: e B}), we have
$$
e^b\Dorfman{\frkd+\mu,\frkr+\nu}_\Theta=\Dorfman{\frkd+\mu,\frkr+\nu},\quad\forall~\frkd+\mu,~\frkr+\nu\in\Gamma(\omni).
$$
Furthermore,  $b$ being skew-symmetric implies  that $e^b$ preserves
the standard pairing (\ref{standard pair}). Therefore, the
transformation $e^b$ is an isomorphism. \qed

\section{Exact $E$-Courant algebroids}
\begin{defi}An $E$-Courant algebroid
$(\huaK,\ppairingE{\cdot,\cdot},\CDorfman{\cdot,\cdot},\rho)$ is
said to be exact if
the following sequence is exact:
\begin{equation}\label{Seq:ExactCourant}
\xymatrix@C=0.5cm{0 \ar[r] & \jet E  \ar[rr]^{\rhowx} &&
                \huaK  \ar[rr]^{\rho} && \dev E \ar[r]  & 0.
                }
\end{equation}
\end{defi}

 \vspace{1mm}
Obviously,   omni-Lie algebroids are   exact. In \cite{PW},
 it is shown that any exact Courant algebroid structure on $TM\oplus T^*M$
 is a
 twist of the standard one  by a closed 3-form.  An important ingredient in the proof of this fact
is that any exact Courant algebroid has an {\bf isotropic
splitting}, i.e. both $TM$ and $T^*M$ are isotropic subbundles.
Unfortunately, this fact is no longer true for an exact $E$-Courant
algebroid when $\mathrm{rank}(E)\geq 2$. Therefore, we shall need
the language of Leibniz algebra cohomologies.

 Recall that a representation of the Leibniz algebra
$(\frkg,[\cdot,\cdot])$ is an $R$-module $V$ equipped with,
respectively, left and right actions   of $\frkg$ on $V$,
$$[\cdot,\cdot]:\frkg\otimes V\longrightarrow V,\quad [\cdot,\cdot]:V\otimes\frkg \longrightarrow V,$$
such that
$$l_{[g_1,g_2]}=[l_{g_1},l_{g_2}],\quad r_{[g_1,g_2]}=[l_{g_1},r_{g_2}],\quad r_{g_2}\circ l_{g_1}=-r_{g_2}\circ r_{g_1},$$
where $l_{g_1}g=[g_1,g]$ and $r_{g_1}g=[g,g_1]$. The Leibniz
cohomology of $\frkg$ with coefficients in $V$ is the homology of
the cochain complex $C^k(\frkg,V)=\Hom_R(\otimes^k\frkg,A),
(k\geq0)$ with the coboundary operator
$\partial:C^k(\frkg,V)\longrightarrow C^{k+1}(\frkg,V)$  defined by
\begin{eqnarray*}
\partial c^k(g_1,\cdots,g_{k+1})&=&\sum_{i=1}^k(-1)^{i+1}g_i((c^k(g_1,\cdots,\widehat{g_i},\cdots,g_{k+1}))
+(-1)^{k+1}(c^k(g_1,\cdots,g_k))g_{k+1}\\
&&+\sum_{1\leq i<j\leq
k+1}(-1)^ic^k(g_1,\cdots,\widehat{g_i},\cdots,g_{j-1},[g_i,g_j],g_{j+1},\cdots,g_{k+1}).
\end{eqnarray*}
The fact that $\partial\circ\partial=0$ is proved in \cite{Loday and
Pirashvili}.

For the omni-Lie algebroid $\omni$,  $\Gamma(\omni)$ is a Leibniz
algebra \cite{clomni} and $\Gamma(\dev E)$ is a Lie algebra. So,
there are two actions of the Leibniz algebra $\Gamma(\dev E)$ on
$\Gamma(\jet E)$, respectively, defined by
\begin{equation}\label{eqn:Leibniz action}
\Dorfman{\frkd,\nu}=\Liederivative_{\frkd}\nu,\quad
\Dorfman{\mu,\frkr}=-\Liederivative_{\frkr}\mu +
\jetd\conpairing{\mu,\frkr}_E\,,\quad \forall~\frkd,~\frkr\in
\Gamma(\dev E), ~\mu,~\nu\in \Gamma(\jet E).
\end{equation}
For any $b\in\Gamma(\Hom(\dev E,\jet E))$, $\partial b=0$ is
equivalently saying that
$$
\Liederivative_{\frkd}b(\frkr)-\Liederivative_{\frkr}b(\frkd) +
\jetd\conpairing{b(\frkd),\frkr}_E-b([\frkd,\frkr]_\dev)=0.
$$
For any $\Theta\in \Gamma(\Hom(\otimes^2\dev E,\jet E))$, $\partial
\Theta=0$ is equivalently saying that
\begin{equation}\label{eqn:two cocycle}
\Liederivative_\frkd\Theta(\frkr,\frkt)-\Liederivative_\frkr\Theta(\frkd,\frkt)+\Liederivative_\frkt\Theta(\frkd,\frkr)
-\jetd\langle\frkt,\Theta(\frkd,\frkr)\rangle_E+\Theta(\frkd,[\frkr,\frkt]_\dev)-\Theta(\frkr,[\frkd,\frkt]_\dev)
-\Theta([\frkd,\frkr]_\dev,\frkt)=0.
\end{equation}
In the meantime, we treat the     Lie algebra  $\Gamma(\dev E)$ as a
Leibniz algebra and define left and right actions of $\Gamma(E)$ on
$E$, respectively, by
\begin{equation}\label{eqn:Leibniz action 2}
[\frkd,u]=\frkd u,\quad [v,\frkr]=-\frkr v\quad
\forall~\frkd,~\frkr\in \Gamma(\dev E), ~u,~v\in\Gamma(E).
\end{equation}
Note that any 2-chain $\Theta\in \Gamma(\Hom(\otimes^2\dev E,\jet
E))$ can be considered as a 3-chain $\widehat{\Theta}\in
\Gamma(\Hom(\otimes^3\dev E, E))$:
$$
\widehat{\Theta}(\frkd\otimes\frkr\otimes\frkt)=\conpairing{\Theta(\frkd\otimes\frkr),\frkt}_E.
$$
The following lemma can be easily verified.
\begin{lem}
The above notation being maintained,   $\Theta\in
\Gamma(\Hom(\otimes^2\dev E,\jet E))$ is closed in the cohomology of
the Leibniz algebra $\Gamma(\dev E)$ with  coefficients in
$\Gamma(\jet E)$ if and only if $\widehat{\Theta}\in
\Gamma(\Hom(\otimes^3\dev E,E))$ is
  closed in the cohomology of  the Leibniz algebra $\Gamma(\dev
E)$ with  coefficients in $\Gamma(E)$.

\end{lem}

\begin{defi}
   Given a  symmetric bundle map $\omega:\dev E\otimes\dev
E\longrightarrow E$  and $\Theta\in \Gamma(\Hom(\otimes^2\dev E,\jet
E))$, the  pair $(\omega,\Theta)$ is called an {\bf{admissible
pair}} of the omni-Lie algebroid $\omni$ if the following conditions
are satisfied:
\begin{itemize}
\item[1) ]$\Theta $ is a 2-cocycle of the Leibniz cohomology of
$\Gamma(\dev E)$ with coefficients in $\Gamma(\jet E)$;
\item[2) ]for any $\frkd\in\Gamma(\dev E)$,
$\Theta(\frkd\otimes\frkd)=\jetd(\omega(\frkd\otimes\frkd))$;
\item[3) ]for any $\frkd,~\frkr\in\Gamma(\dev E)$,
$\half\frkd\omega(\frkr\otimes\frkr)=
\ppairingE{\Theta(\frkd\otimes\frkr),\frkr}+\omega([\frkd,\frkr]_\dev\otimes\frkr).$
\end{itemize}
Two admissible  pairs $(\omega,\Theta)$ and
$(\widetilde{\omega},\widetilde{\Theta})$ are said to be equivalent
 if there is some
$b\in\Gamma(\Hom(\dev E,\jet E))$ such that
\begin{itemize}
\item[1)]for any $\frkd,~\frkr\in\Gamma(\dev E)$,
$\widetilde{\omega}(\frkd\otimes\frkr)=\omega(\frkd\otimes\frkr)
+\frac{1}{2}(\conpairing{b(\frkd),\frkr}_E+\conpairing{b(\frkr),\frkd}_E)$;
\item[2)]$\widetilde{\Theta}=\Theta+\partial b$.
\end{itemize}
\end{defi}
For every $b\in\Gamma(\Hom(\dev E,\jet E))$, $\partial b$ is a
2-cocycle and we can define a symmetric bundle map $\omega_b:\dev
E\otimes\dev E\longrightarrow E$   by
$$
\omega_b(\frkd\otimes\frkr)=\frac{1}{2}(\conpairing{b(\frkd),\frkr}_E+\conpairing{b(\frkr),\frkd}_E).
$$
  It is straightforward to verify that the
pair $(\omega_b,\partial b)$ is an admissible  pair and it is
equivalent to the admissible  pair $(0,0)$.

\begin{thm}\label{Thm:twirsted}
There is a one-to-one correspondence between   isomorphic classes of
exact $E$-Courant algebroids and equivalence classes of admissible
pairs of the omni-Lie algebroid $\omni= \dev E\oplus\jet E$.

More precisely, for any exact $E$-Courant algebroid
$(\huaK,\ppairingE{\cdot,\cdot},\CDorfman{\cdot,\cdot},\rho)$, one
may identify $\huaK=\dev E\oplus\jet E$ and  take $\rho$ as the
projection to $\dev E$, and then there is a corresponding admissible
pair $(\omega,\Theta)$ such that
\begin{eqnarray}
\label{Eqn:Twistedpair}\ppairingE{\frkd+\mu,\frkr+\nu}&=&\frac{1}{2}(\conpairing{\frkd,\nu}_E
+\conpairing{\frkr,\mu}_E)+\omega(\frkd\otimes\frkr),\quad\forall~ \frkd+\mu,\frkr+\nu\in\Gamma(\dev E\oplus\jet E),\\
\label{Eqn:TwistedTheta} \CDorfman{\frkd+\mu,\frkr+\nu}
&=&\Dorfman{\frkd+\mu,\frkr+\nu}_\Theta=\Dorfman{\frkd+\mu,\frkr+\nu}+
\Theta(\frkd\otimes \frkr),
\end{eqnarray}
where $\Dorfman{\cdot,\cdot}$ is the standard bracket
(\ref{standard bracket}).

Conversely, for any admissible  pair $(\omega,\Theta)$, $(\dev
E\oplus\jet E,\ppairingE{\cdot,\cdot},\CDorfman{\cdot,\cdot},\rho)$
is an exact $E$-Courant algebroid, where $\ppairingE{\cdot,\cdot}$
and $\CDorfman{\cdot,\cdot}$ are given as above. Moreover, two exact
$E$-Courant algebroids are isomorphic if and only if their
corresponding  admissible pairs are equivalent.
\end{thm}
\pf We split the proof into four steps. In Step 1, we prove that the
$E$-valued pairing of the exact $E$-Courant algebroid $\huaK$ is
given by (\ref{Eqn:Twistedpair}). In Step 2, we prove that the
bracket of $\huaK$ is given by (\ref{Eqn:TwistedTheta}).   In Step
3, we prove that if we choose different splitting, we obtain
equivalent admissible pairs. In Step 4, we give the proof of the
reverse statement.

{\bf{Step 1.}} By choosing a splitting $s:\dev E\longrightarrow
\huaK$ of the exact sequence (\ref{Seq:ExactCourant}), we have
   $\huaK\cong\omni=\dev
E\oplus \jet E$  and $\rho$ is then  the projection from $\omni$ to
$\dev E$.
By Properties (EC-4) and (EC-5), for all $\mu,~\nu\in\jet E$, we
have
\begin{equation}\label{eqn:pair of jet}
\ppairingE{\rhowx\mu,\rhowx\nu}=\frac{1}{2}\conpairing{\mu,\rho\circ\rhowx\nu}_E=0,
\end{equation}
which implies that  $\rhowx\jet E$ is isotropic under the pairing
$\ppairingE{\cdot,\cdot}$. So if we transfer the pairing
$\ppairingE{\cdot,\cdot}$ on $\huaK$ to a pairing
$\ppairingE{\cdot,\cdot}$ on $\dev E\oplus \jet E$, $\jet E $ is
isotropic. For any $\frkd\in\dev E$, we have
\begin{equation}\label{eqn:pair of jet and dev}
\ppairingE{\frkd,\nu}=\ppairingE{s(\frkd),\rhowx(\nu)}
=\frac{1}{2}\conpairing{\frkd,\nu}_E.
\end{equation}
Furthermore, for all $\frkd,\frkr\in \dev E$, we have
\begin{equation}\label{eqn:pair of dev}
\ppairingE{\frkd,\frkr}=\ppairingE{s(\frkd),s(\frkr)}\triangleq
\omega(\frkd\otimes \frkr),
\end{equation}
where $\omega:\dev E\otimes\dev E\longrightarrow E$ is a symmetric
bundle map. By $(\ref{eqn:pair of jet}),~(\ref{eqn:pair of jet and
dev})$ and (\ref{eqn:pair of dev}), it follows that the pairing  is
given by (\ref{Eqn:Twistedpair}).

{\bf{Step 2.}}  ~For any $\frkd$, $\frkr\in\Gamma(\dev E)$, by
Properties (EC-1) and (EC-2), we are able to write
\begin{equation}\label{Eqt:tp21}
\CDorfman{\frkd,\frkr}=[\frkd,\frkr]_\dev+\Theta({\frkd}\otimes
{\frkr}),
\end{equation}
where $\Theta$ is an $\mathbb R$-linear mapping $\Gamma( \dev E)
\otimes\Gamma (\dev E)\lon \Gamma (\jet E)$. By (\ref{eqn:bracket X
fY}), we know that $\Theta$ is also $\CWM$-linear and hence
$\Theta\in \Gamma (\Hom(\otimes^2\dev E, \jet E))$.

Again by Property (EC-1), there is a bi-linear map
$$
\Delta: \Gamma (\dev E)\times \Gamma( \jet E)\lon \Gamma (\jet E),
$$
such that
$$
\CDorfman{\frkd,\mu}=\Delta(\frkd,\mu),\quad\forall ~~
\frkd\in\Gamma (\dev E),~\mu\in\Gamma (\jet E).
$$
By Property (EC-3) and $\jet E$ being isotropic, we have
$$
\frkd\ppairingE{\frkr,\mu}  =\ppairingE{\CDorfman{\frkd,\frkr} ,\mu}
+\ppairingE{\frkr,\CDorfman{ \frkd,\mu} }
=\ppairingE{[\frkd,\frkr]_\dev,\mu}+\ppairingE{\frkr,\Delta(\frkd,\mu)},
$$
which implies that $\Delta(\frkd,\mu)=\Liederivative_{\frkd}\mu$,
i.e.
\begin{equation}\label{Eqt:tp22}
\CDorfman{\frkd,\mu}=\Liederivative_{\frkd}\mu,\quad\forall ~~
\frkd\in\Gamma (\dev E),~\mu\in\Gamma (\jet E).
\end{equation}
Furthermore, we have
\begin{eqnarray*}
\jetd \conpairing{\frkd,\mu}_E&=&\jetd
\ppairingE{\frkd+\mu,\frkd+\mu}-\jetd\ppairingE{\frkd,\frkd}=\CDorfman{\frkd+\mu,\frkd+\mu}-\jetd\ppairingE{\frkd,\frkd}
\\&=&\CDorfman{\frkd,\mu}+\CDorfman{\mu,\frkd}
=\Liederivative_{\frkd}\mu+\CDorfman{\mu,\frkd}.
\end{eqnarray*}
Therefore,
\begin{equation}\label{Eqt:tp23}
\CDorfman{\mu,\frkd} =\jetd
\conpairing{\frkd,\mu}_E-\Liederivative_{\frkd}\mu,\quad\forall ~
~\mu\in\Gamma (\jet E),\frkd\in\Gamma (\dev E).
\end{equation}
 Again by Property (EC-3), we have
$$
\ppairingE{\CDorfman{\mu,\nu},\frkd}+\ppairingE{\nu,\CDorfman{\mu,\frkd}}=0,\quad
\forall~ \mu,\nu\in\Gamma(\jet E),\frkd\in \Gamma(\dev E),
$$
which implies
\begin{equation}\label{eqt:1111}
\CDorfman{\mu,\nu}=0,\quad \forall ~\mu,~\nu\in\Gamma(\jet E).
\end{equation}
By  (\ref{Eqt:tp21}), (\ref{Eqt:tp22}), (\ref{eqt:1111}) and
(\ref{Eqt:tp23}), we get
\begin{equation}\label{Eqt:DorfmanUpsilon}
\CDorfman{\frkd+\mu,\frkr+\nu} =\Dorfman{\frkd+\mu,\frkr+\nu}+
\Theta({\frkd}\otimes{\frkr}),\quad\forall~
\frkd+\mu,\frkr+\nu\in\Gamma(\omni).
\end{equation}
Since the bracket $\CDorfman{\cdot,\cdot}$ satisfies the Leibniz
rule, we have
$$
\Liederivative_\frkd\Theta(\frkr,\frkt)-\Liederivative_\frkr\Theta(\frkd,\frkt)+\Liederivative_\frkt\Theta(\frkd,\frkr)
-\jetd\langle\frkt,\Theta(\frkd,\frkr)\rangle_E+\Theta(\frkd,[\frkr,\frkt]_\dev)-\Theta(\frkr,[\frkd,\frkt]_\dev)
-\Theta([\frkd,\frkr]_\dev,\frkt)=0,
$$
which implies that $\Theta $ is a 2-cocycle in the Leibniz
cohomology of $\Gamma(\dev E)$ with coefficients in $\Gamma(\jet E)$
and the two actions are given in (\ref{eqn:Leibniz action}). Since
the pair $(\omega,\Theta)$ comes from the $E$-Courant algebroid
$\huaK$, it is straightforward to verify  that it is an admissible
pair.

{\bf{Step 3.}}~Suppose that we have two sections $s_1,~s_2:\dev
E\longrightarrow\huaK$, then  $\rho(s_1-s_2)=0$. Take
$b=s_1-s_2:\dev E\longrightarrow\jet E$. We then have
$s_2(\frkd)=\frkd+b(\frkd)$. By straightforward computations, we get
\begin{eqnarray}
\nonumber\Dorfman{\frkd+b(\frkd),\frkr+b(\frkr)}_\Theta&=&[\frkd,\frkr]_\dev
+\Liederivative_\frkd
b(\frkr)-i_\frkr\partial b(\frkd)+\Theta(\frkd\otimes\frkr)\\
\label{eqn:splitting}&=&[\frkd,\frkr]_\dev+b([\frkd,\frkr]_\dev)+(\Theta+\partial
b)(\frkd\otimes\frkr).
\end{eqnarray}
  If we denote the new pairing by $\widetilde{\omega}$,
for any $\frkd,~\frkr\in\Gamma(\dev E)$, then
$$
\widetilde{\omega}(\frkd,\frkr)=\omega(\frkd,\frkr)+\ppairingE{b(\frkd),\frkr}+\ppairingE{b(\frkr),\frkd}.
$$
Therefore, if we choose different splitting, we obtain equivalent
admissible  pairs.

{\bf{Step 4.}}~ Conversely, for any admissible  pair
$(\omega,\Theta)$, on $\dev E\oplus\jet E$, we define the pairing
$\ppairingE{\cdot,\cdot}$ and the bracket
$\Dorfman{\cdot,\cdot}_\Theta$ by (\ref{Eqn:Twistedpair}) and
(\ref{Eqn:TwistedTheta}). It is straightforward to see that $(\dev
E\oplus\jet
E,\ppairingE{\cdot,\cdot},\Dorfman{\cdot,\cdot}_\Theta,\rho)$ is an
$E$-Courant algebroid. If we choose different representative element
$(\widetilde{\omega},\widetilde{\Theta})$, assume that
$\widetilde{\Theta}=\Theta+\partial b$ for some
$b\in\Gamma(\Hom(\dev E,\jet E))$ and the corresponding $E$-Courant
algebroid is $(\dev E\oplus\jet
E,\ppairingE{\cdot,\cdot}^\prime,\Dorfman{\cdot,\cdot}_{\Theta+\partial
b},\rho)$. By some computations similar to (\ref{eqn:splitting}),
for any $\frkd+\mu,~\frkr+\nu\in\Gamma(\dev E\oplus\jet E)$, we have
$$
e^b\Dorfman{\frkd+\mu,\frkr+\nu}_{\Theta+\partial
b}=\Dorfman{e^b(\frkd+\mu),e^b(\frkr+\nu)}_{\Theta}.
$$
It is also obvious  that
$$
\ppairingE{e^b(\frkd+\mu),e^b(\frkr+\nu)}=\ppairingE{\frkd+\mu,\frkr+\nu}^\prime.
$$
Therefore, the transformation $e^b$ is the isomorphism from $(\dev
E\oplus\jet
E,\ppairingE{\cdot,\cdot}^\prime,\Dorfman{\cdot,\cdot}_{\Theta+\partial
b},\rho)$ to $(\dev E\oplus\jet
E,\ppairingE{\cdot,\cdot},\Dorfman{\cdot,\cdot}_\Theta,\rho)$. \qed

\vspace{3mm}

\begin{rmk}
The extreme case that the induced symmetric bundle map $\omega:\dev
E\otimes \dev E\longrightarrow E$ is zero, i.e. the splitting is
isotropic,  already  has been studied at the end of Section 3, which
is in fact the twisted omni-Lie algebroid.
\end{rmk}
In some special cases, we can define an isotropic splitting as
follows.
\begin{pro}Under the circumstances above, if the induced $E$-valued pairing
$\omega:\dev E\otimes\dev E\longrightarrow E$ on $\dev E$ satisfies
 $\Img(\omega_\natural)\subset \jet E$,  then there is an isotropic
splitting  $s(\frkd)=-\frac{1}{2}\omega_\natural(\frkd)$.
In particular, if $E$ is a line bundle, there  always exist
isotropic splittings.
\end{pro}
\pf    By definition,  we have
$$
\ppairingE{\frkd+s(\frkd),\frkr+s(\frkr)}=\ppairingE{\frkd,\frkr}-\omega(\frkd,\frkr)=0,
$$
for all $\frkd,~\frkr\in\dev E$ and $\mu,~\nu\in\jet E$. Thus,
$s(\dev E)$ is isotropic and we proved the first claim.

Moreover, for any $\Phi\in\gl(E)$,
$$
\Img(\omega_\natural)\subset \jet
E\Longleftrightarrow\omega(\frkd,\Phi)=\Phi(\omega(\frkd,\Id_E)).
$$
 If E is a line bundle, the conclusion follows, because $\gl(E)$
is then a trivial line bundle. \qed

\section{$E$-Lie bialgebroids}
In this section we introduce the notion of an $E$-Lie bialgebroid,
whose double turns out to be an $E$-Courant algebroid. Conversely,
any $E$-Courant algebroid which is the direct sum of two transverse
Dirac structures provides an $E$-Lie bialgebroid. Similar to the
fact that the base manifold of a Lie bialgebroid is a Poisson
manifold, for an $E$-Lie bialgebroid, the underlining vector bundle
$E$ is a Lie algebroid (if $\mathrm{rank} E\geq2$), or a local Lie
algebra (if $\mathrm{rank} E=1$).

In the sequel,  notations introduced in Section 3 are needed.

\begin{defi}\label{Def:Ebialgebroid} An $E$-dual pair $((A,\rho_A);(B,\rho_B))$ is called  an $E$-Lie bialgebroid if for
all $X,~Y\in \Gamma(A)$, $u,~v\in \Gamma(E)$, the following
conditions are satisfied:
\begin{itemize}
\item[(1)]$\dB[X,Y]=\Liederivative_{X }(\dB Y)-\Liederivative_{Y}(\dB
X )$,
\item[(2)] $\Liederivative_{\dA u}X=-\Liederivative_{\dB u}X$,
\item[(3)] $\conpairing{\dB u, \dA u}_E=0.$
\end{itemize}
\end{defi}When there is no confusion,   we simply denote such an $E$-Lie
bialgebroid by $(A,B)$.
\begin{rmk}
Condition $(3)$ is equivalent to  $\rho_B\circ \dA=-\rho_A\circ
\dB$.
\end{rmk}
Let us give some examples. Recall   the properties of an omni-Lie
algebroid $ \omni=\dev{E}\oplus \jet{E}$. The pair
 $(\dev E,\jet E)$ is certainly an example of $E$-Lie bialgebroids, where $\rho_{\jet E}=0$
and $\rho_{\dev E}$ is the identity map. For a Lie bialgebroid
$(A,A^*)$, it is an $E$-Lie bialgebroid, where $E$ is the trivial
line bundle $M\times \mathbb R$. The representations $\rho_A$ and
$\rho_{A^*}$ are, respectively, the anchors of $A$ and $A^*$. For a
Lie algebroid $A$ and a representation $\rho_A:A\longrightarrow\dev
E$, $(A,A^*\otimes E) $
 is an $E$-Lie bialgebroid, where $\rho_{A^*\otimes E}=0$ (Example \ref{ex: A A dual tensor E}).

\begin{pro}
A generalized Lie bialgebroid $((A,\phi_0),(A^*,X_0))$ is an $E$-Lie
bialgebroid, where $E $ is the trivial line bundle $M\times \mathbb
R$.
\end{pro}
\pf For any $X\in \Gamma(A),~\xi\in\Gamma( A^*)$, the
representations $\rho_A$ and $\rho_{A^*}$ are given by
$$\rho_A(X)=a(X)+\phi_0(X),\quad \rho_{A^*}(\xi)=a_*(\xi)+X_0(\xi),$$
where $a$ and $a_*$ are, respectively, the anchors of $A$ and $A^*$.
Evidently, we have $\dM^{A^*}=d_{*X_0}$. Furthermore, by definition,
we have
\begin{equation}\label{eqn:con 1}
[X,d_{*X_0}Y]_{\phi_0}=[X,d_{*X_0}Y]-\langle\phi_0,X \rangle
d_{*X_0}Y=\Liederivative_Xd_{*X_0}Y,\quad \forall~X,Y\in\Gamma(A).
\end{equation} Since
$((A,\phi_0),(A^*,X_0))$ is  a generalized Lie bialgebroid, we have
 (\ref{eqn:gene Lie bi 1}).  By (\ref{eqn:con 1}),
Condition (1) of Definition \ref{Def:Ebialgebroid} holds, i.e.
\begin{equation}\label{eqn:con 11}
\dM^{A^*}[X,Y]=\Liederivative_{X }(\dM^{A^*}
Y)-\Liederivative_{Y}(\dM^{A^*} X ).
\end{equation}
To prove  Condition (2) of Definition \ref{Def:Ebialgebroid}, one
substitutes $Y$ by $fY$ in (\ref{eqn:con 11}), where $f\in
C^\infty(M)$, and gets
$$
\Liederivative_{\dM^Af}X=-\Liederivative_{\dM^{A^*}f}X+f(\Liederivative_{*\phi_0}X+\Liederivative_{X_0}X).
$$
By   (\ref{eqn:gene Lie bi 2}), we have
\begin{equation}\label{eqn:con 2}
\Liederivative_{\dM^Af}X=-\Liederivative_{\dM^{A^*}f}X,
\end{equation}
which is exactly Condition (2) of Definition \ref{Def:Ebialgebroid}.
Finally, substituting  $X$ by $fX$ in (\ref{eqn:con 2}), we get
$$
\rho_{A^*}\dM^A(f)=-\rho_{A}\dM^{A^*}(f)+(a(X_0)+a_*(\phi_0))(X).
$$
By (\ref{eqn:gene Lie bi 2}) again, we have
$\rho_{A^*}\dM^A(f)=-\rho_{A}\dM^{A^*}(f)$, which is exactly
Condition (3) of Definition \ref{Def:Ebialgebroid}. Therefore, a
generalized Lie bialgebroid $((A,\phi_0),(A^*,X_0))$ is truly an
$E$-Lie bialgebroid. \qed \vspace{3mm}

Let $(A,[\cdot,\cdot],a)$ be a Lie algebroid and $\rho_A: ~A\lon
\dev E$ a $B$-invariant representation, where $B$ is an $E$-dual
bundle of $A$. For any $u\in \Gamma(E),~X\in\Gamma(A)$ and
$X^k\in\Gamma(\Hom(\wedge^k B,E)_A)$, we  define their Schouten
brackets by
$$[u,X^k]=[X^k,u]=(-1)^{k+1}i_{\dA u}X^k,\quad [X,X^k]=-[X^k,X]=\Liederivative_XX^k.$$ The Schouten bracket $[H,K]\in \Gamma(\Hom(\wedge^3B,E)_A)$ of
$H,~K\in \Gamma(\Hom(\wedge^2B,E)_A)$ is  defined by
\begin{equation}
[H,K](\xi_1,\xi_2,\xi_3)=\langle
\Liederivative_{K\xi_1}\xi_2,H\xi_3\rangle_E+\langle
\Liederivative_{H\xi_1}\xi_2,K\xi_3\rangle_E+c.p.,\quad
\forall~\xi_i\in \Gamma(B).
\end{equation}
For any $\Lambda\in\Gamma(\Hom(\wedge^2B,E)_A)$, we introduce a
bracket $[\cdot,\cdot]_\Lambda$ on $\Gamma(B)$:
\begin{equation}\label{def:bracket}
[\xi,\eta]_\Lambda=\Liederivative_{\Lambda\xi}\eta-\Liederivative_{\Lambda\eta}\xi-\dA
(\Lambda(\xi,\eta)),\quad\forall~\xi,~\eta\in\Gamma(B).
\end{equation}
By straightforward computations, we have the following formula:
\begin{equation}\label{eqn: H}
[H\xi,H\eta]=H[\xi,\eta]_H+\frac{1}{2}[H,H](\xi,\eta),
\end{equation}
which implies that $[H,H]\in \Gamma(\Hom(\wedge^3B,E)_A)$. Moreover,
replacing $H$ by $H+K$, we know that
$[H,K]\in\Gamma(\Hom(\wedge^3B,E)_A)$. \vspace{3mm}

Given  $\Lambda\in\Gamma(\Hom(\wedge^2B,E)_A)$, let
$\Lambda_\natural:B\longrightarrow A$ and
$[\Lambda,\Lambda]_\natural: \wedge^2B \longrightarrow A$ be the
induced bundle maps defined by (\ref{map}). Let us denote
$$ a_B = a\circ\Lambda_\natural:B\longrightarrow TM \quad \mbox{and} \quad
 \rho_B = \rho_A\circ\Lambda_\natural:B\longrightarrow\dev E.$$

\begin{pro}\label{pro:B}
Under the circumstances above,  $(B,[\cdot,\cdot]_\Lambda,a_B)$ is a
Lie algebroid together with an $A$-invariant representation $\rho_B$
if and only if the following two conditions are satisfied:
 \begin{itemize}
\item[(1)]$\rho_A\circ[\Lambda,\Lambda]_\natural=0;$
\item[(2)]$\Liederivative_X[\Lambda,\Lambda]=0,\quad\forall~X\in\Gamma(A).$
\end{itemize}
\end{pro}
\pf By (\ref{eqn: H}), for any $\xi,~\eta\in \Gamma(B)$, we have
\begin{equation}\label{eqn:representation}
\rho_B([\xi,\eta]_\Lambda)=\rho_A\circ\Lambda_\natural([\xi,\eta]_\Lambda)=[\rho_B\xi,\rho_B\eta]
-\frac{1}{2}\rho_A\circ[\Lambda,\Lambda]_\natural(\xi,\eta).
\end{equation}
Therefore, $\rho_B$ is a homomorphism if and only if
$\rho_A\circ[\Lambda,\Lambda]_\natural=0.$ It is simple to see that
for any $f\in C^\infty(M)$,
$$
[\xi,f\eta]_\Lambda=a_B(\xi)(f)\eta+f[\xi,\eta]_\Lambda.
$$
Furthermore, by (\ref{eqn: H}), we have
$$
a_B([\xi,\eta]_\Lambda)=[a_B\xi,a_B\eta]-\frac{1}{2}\jd\circ\rho_A\circ[\Lambda,\Lambda]_\natural(\xi,\eta).
$$
Therefore, if $\rho_A\circ[\Lambda,\Lambda]_\natural=0,$ $a_B$ is a
homomorphism.

    Let
$$
J(\xi_1,\xi_2,\xi_3)=[[\xi_1,\xi_2]_\Lambda,\xi_3]_\Lambda+c.p.,\quad
\forall~\xi_i\in\Gamma(B).
$$
For any $X\in\Gamma(A)$,  under the condition that
$\rho_A\circ[\Lambda,\Lambda]_\natural=0$ and by some similar
calculations did in \cite{L X}, one is able to get
\begin{eqnarray*}
\conpairing{J(\xi_1,\xi_2,\xi_3),X}_E&=&-\{\frac{1}{2}\conpairing{[X,\Lambda_\natural(\xi_1,\xi_2)],\xi_3}_E+c.p.\}
+\rho_A(X)([\Lambda,\Lambda](\xi_1,\xi_2,\xi_3))\\
&=&-\frac{1}{2}(\Liederivative_X[\Lambda,\Lambda])(\xi_1,\xi_2,\xi_3).
\end{eqnarray*}
Thus, under the condition that
$\rho_A\circ[\Lambda,\Lambda]_\natural=0$, the bracket
$[\cdot,\cdot]_\Lambda$ defined by (\ref{def:bracket}) satisfies the
Jacobi identity if and only if $\Liederivative_X[\Lambda,\Lambda]=0$
for all $~X\in\Gamma(A).$

Finally, we show that the representation $\rho_B$ of the Lie
algebroid $(B,[\cdot,\cdot]_\Lambda,a_B)$ is $A$-invariant. In fact,
by straightforward computations, we have
\begin{equation}\label{eqn:d B}
\dB u=-i_{\dA u}\Lambda,\quad\dB
X=-\Liederivative_X\Lambda,\quad\forall
~u\in\Gamma(E),~X\in\Gamma(A),
\end{equation}
as required. \qed

\begin{thm}
Under the conditions of Proposition \ref{pro:B}, $(A,B)$ is an
$E$-Lie bialgebroid.
\end{thm}\pf By (\ref{eqn:d B}), it is straightforward to check   the
compatibility conditions of an $E$-Lie bialgebroid. We omit the
details. \qed\vspace{3mm}

The notion of local Lie algebra was introduced by Kirillov in
\cite{KirillovLocal}, which is a vector bundle $E$ whose section
space $\Gamma(E)$ has an $\Real$-Lie algebra structure
$[\cdot,\cdot]_E$ with the local property,
$\mathrm{supp}[u,v]\subset \mathrm{supp}{u}\cap
  \mathrm{supp}{v}$,  for all $u, v \in \Gamma(E)$.
In particular, $M$ is called a Jacobi manifold  if the trivial
bundle $M\times\Real $ is a local Lie algebra, which is equivalent
 to that there is a pair $(\Lambda, X)$, where $\Lambda$ is a
bi-vector field and $X$ is a vector field on $M$ such that
$[\Lambda,\Lambda]=2 X\wedge \Lambda$ and $[\Lambda,X]=0$. A Jacobi
structure reduces to a Poisson structure if $X=0$. Similar to the
fact that the cotangent bundle of a Poisson manifold is a Lie
algebroid, $T^*M\oplus \mathbb R$ enjoys a Lie algebroid structure
for every Jacobi manifold. A local Lie algebra is not a Lie
algebroid since there is no anchor map. For example, the trivial
bundle $M\times \mathbb R$ with the Poisson bracket is only a local
Lie algebra for every Poisson manifold $M$.

It is well known that a Lie bialgebroid $(A,A^*)$  gives a Poisson
structure on the base manifold $M$ (\cite{MackenzieX:1994}). If
$(A,A^*)$ is a generalized Lie bialgebroid, there is an induced
Jacobi structure on the base manifold $M$ (\cite{generalize Lie
bialgebroid}). In the situation of an $E$-Lie bialgebroid $(A,B)$ as
in Definition \ref{Def:Ebialgebroid}, we introduce a bracket
$[\cdot,\cdot]_E$ on $\Gamma(E)$ as follows
\begin{equation}\label{Eqt:Ebracket}
[u,v]_E\defbe \conpairing{\dA u, \dB v}_E ~(~= \rho_B(\dA u)
v~),\quad\forall ~ u,~v\in \Gamma(E).
\end{equation}

\begin{thm}\label{Thm:ElocalLie}Let $(A,B)$ be an $E$-Lie bialgebroid,
\begin{itemize}
\item[1)] If $\mathrm{rank}(E)\geq 2$,  $(E,[\cdot,\cdot]_E)$ is a Lie
algebroid with the anchor $\jd \circ\rho_B\circ
\rhoAWuXing\circ\jetd$.
\item[2)] If $\mathrm{rank}E=1$,  $(E,[\cdot,\cdot]_E)$ is a local Lie algebra.
\end{itemize}
\end{thm}
\pf By Condition (3) of Definition \ref{Def:Ebialgebroid},  the
bracket defined by (\ref{Eqt:Ebracket}) is skew-symmetric. To check
the Jacobi identity, for all $u,~v,~w\in \Gamma(E)$, we have
\begin{eqnarray*}
[u,[v,w]_E]_E&=&\conpairing{\dA u, \dB\conpairing{\dA v, \dB
w}_E}_E=\conpairing{\dA u, \dB (i_{\dA v}\dB w)+ i_{\dA v}\dB \dB w}_E\\
&=&-\conpairing{\dA u, \Liederivative_{\dA v} \dB w}_E=-
\conpairing{\dA u, [\dB v,\dB w]}_E=-\rho_A([\dB v,\dB w])u\\
&=&-\rho_A(\dB v)\rho_A(\dB w)u+\rho_A(\dB w)\rho_A(\dB
v)u\\&=&-[[u,w]_E,v]_E+[[u,v]_E,w]_E \,.
\end{eqnarray*}
Moreover, we have an obvious expression:
\begin{eqnarray}\nonumber
[u,fv]_E&=&\conpairing{\dA u, f\dB v+\rhoAWuXing(\dM f\otimes
v)}_E\\\nonumber &=&f[u,v]_E+\conpairing{\rho_B(\dA u),\dM f\otimes
v}_E\\\label{Eqt:ufvE} &=& \nonumber f[u,v]_E+(\jd \circ\rho_B\circ
\rhoAWuXing(\jetd u))(f) v.
\end{eqnarray}
The map $\rho_B\circ \rhoAWuXing: \quad~\jet E\lon \dev E$  is
skew-symmetric.  In  \cite{clomni}, it is shown that
$\mathrm{rank}(E)\geq 2$ implies that $\jd \circ\rho_B\circ
\rhoAWuXing\circ\jetd:E\longrightarrow TM$ is a bundle map. In this
case, $(E,[\cdot,\cdot]_E)$ is a Lie algebroid, whose anchor is $\jd
\circ\rho_B\circ \rhoAWuXing\circ\jetd$. \qed

\section{Manin Triples }

 In this section we develop the theory of Manin triples of $E$-Lie
bialgebroids analogously to that of Lie bialgebroids. It is known
that for a Lie bialgebroid $(A,A^*)$, we can endow $A\oplus A^*$
with a Courant algebroid structure \cite{LWXmani}. For a generalized
Lie bialgebroid $(A,A^*)$, we can endow $A\oplus A^*$ with a
generalized Courant algebroid structure \cite{generalized}.
Similarly, for any $E$-Lie bialgebroid $(A,B)$, we can endow
$A\oplus B$ with an $E$-Courant algebroid structure. In fact, we
define on $A\oplus B$ an $E$-valued pairing
$\ppairingE{\cdot,\cdot}$   by
\begin{equation}\label{eqt:ABpair}
 \ppairingE{X_1+\xi_1,X_2+\xi_2}=\frac{1}{2}(\conpairing{X_1,\xi_2}_E+\conpairing{X_2,\xi_1}_E),
\quad \forall~X_i+\xi_i\in\Gamma(A\oplus B),~i=1,2,
\end{equation}
and a bracket $[\cdot,\cdot]_{A\oplus B}$   by
\begin{eqnarray}\label{bracket1}
\nonumber[X_1+\xi_1,X_2+\xi_2]_{A\oplus
B}&=&[X_1,X_2]+\Liederivative_{\xi_1}X_2-\Liederivative_{\xi_2}X_1+\dB\langle
X_1,\xi_2\rangle_E\\&&+[\xi_1,\xi_2]+\Liederivative_{X_1}\xi_2-\Liederivative_{X_2}\xi_1+\dA\langle
X_2,\xi_1\rangle_E.
\end{eqnarray}
\begin{lem}\label{lem:jacobi}
Let $(A,B)$ be an $E$-Lie bialgebroid. Then, one has
\begin{equation} \label{jacobi}J(e_1,e_2,e_3)=-J_1(e_1,e_2,e_3)+c.p.(e_1,e_2,e_3)-J_2(e_1,e_2,e_3),
\end{equation}
  for all
$~e_i=X_i+\xi_i\in \Gamma(A\oplus B)$. Here
$$J(e_1,e_2,e_3)=\{e_1,\{e_2,e_3\}\}-\{\{e_1,e_2\},e_3\}-\{e_2,\{e_1,e_3\}\}
.$$ The notation c.p. means cyclic permutations.  $J_1(e_1,e_2,e_3)$
and $J_2(e_1,e_2,e_3)$ are, respectively,
\begin{eqnarray}
\nonumber&&J_1(e_1,e_2,e_3)\\
&=&i_{X_3}(\dA[\xi_1,\xi_2]-\Liederivative_{\xi}\dA\xi_2+\Liederivative_{\xi_2}\dA\xi_1)+i_{\xi_3}(\dB[X_1,X_2]
-\Liederivative_{X_1}\dB X_2+\Liederivative_{X_2}\dB
X_1),\\\nonumber
\nonumber&&J_2(e_1,e_2,e_3)\\\nonumber&=&\Liederivative_{\dB\langle\xi_2,X_1\rangle_E}\xi_3+[\dA\langle\xi_2,X_1\rangle_E,\xi_3]
+\Liederivative_{\dB\langle\xi_3,X_1\rangle_E}\xi_2+[\dA\langle\xi_3,X_1\rangle_E,\xi_2]\\
&&
+\Liederivative_{\dB\langle\xi_3,X_2\rangle_E}\xi_1+[\dA\langle\xi_3,X_2\rangle_E,\xi_1].
\end{eqnarray}
\end{lem}
\pf  We need the following formula, $$
i_X\Liederivative_{\xi}\dA\eta=[\xi,\Liederivative_X\eta]-\Liederivative_{\Liederivative_{\xi}X}\eta+
[\dA\langle\eta,X\rangle_E,\xi]+\dA(\rho_B(\xi)\langle\eta,X\rangle_E)-\dA\langle[\xi,\eta],X\rangle_E\,,
$$
for  all $~X\in \Gamma(A),~\xi,~\eta\in \Gamma(B)$. Then, the rest
of the calculations are very similar to those in \cite{LWXmani}.
  We omit the
  details. \qed
\begin{rmk}Equation
(\ref{jacobi}) is different from Theorem 3.1 in \cite{LWXmani} since
the bracket $[\cdot,\cdot]_{A\oplus B}$ given by (\ref{bracket1}) is
not skewsymmetric.
\end{rmk}
\begin{lem}\label{lem:rho}
For any $X\in \Gamma(A),~\xi\in\Gamma(B)$ and $u\in\Gamma(E)$, we
have
\begin{eqnarray*}
\nonumber&&[\rho_B(\xi),\rho_A(X)]_\dev(u)\\
&=&\rho_A(\Liederivative_{\xi}X)(u)-\rho_B(\Liederivative_X\xi)(u)
+\rho_B\circ\dA\langle\xi,X\rangle_E(u)+\langle\Liederivative_{\dB
u}\xi+[\dA u,\xi],X\rangle_E.
\end{eqnarray*}
\end{lem}
\pf
\begin{eqnarray*}
\nonumber&&\rho_B\circ\dA\langle\xi,X\rangle_E(u)\\&=&(\rho_A(\dB
u))\langle\xi,X\rangle_E=\langle\Liederivative_{\dB u}
\xi,X\rangle_E+\langle\xi,[\dB u,X]\rangle_E\\\nonumber
&=&\langle\Liederivative_{\dB u}\xi+[\dA
u,\xi],X\rangle_E+\langle[\xi,\dA
u],X\rangle_E-\langle\xi,\Liederivative_X \dB u\rangle_E\\
&=&[\rho_B(\xi),\rho_A(X)]_\dev(u)-\rho_A(\Liederivative_{\xi}X)(u)+\rho_B(\Liederivative_X\xi)(u)+\langle\Liederivative_{\dB
u}\xi+[\dA u,\xi],X\rangle_E. \qed
\end{eqnarray*}
\begin{thm}
Given an $E$-Lie bialgebroid $(A,B)$ as above, the quadrable
$(A\oplus B, \ppairingE{\cdot,\cdot},[\cdot,\cdot]_{A\oplus
B},\rho_A+\rho_B)$ is an $E$-Courant algebroid, where
$\ppairingE{\cdot,\cdot}$ is given by  (\ref{eqt:ABpair}) and
$[\cdot,\cdot]_{A\oplus B}$  is given by  (\ref{bracket1}).
\end{thm}
\pf
 By Conditions (1), (2) of Definition \ref{Def:Ebialgebroid} and Lemma
 \ref{lem:jacobi},
 $(\Gamma(A\oplus B),[\cdot,\cdot]_{A\oplus B})$ is a Leibniz algebra. For any $X,~Y\in\Gamma(A)$
 and $\xi,~\eta\in\Gamma(B)$, it is obvious that
 $$\rho_A[X,Y]_{A\oplus
B}=[\rho_A(X),\rho_A(Y)]_\dev,\quad \rho_B[\xi,\eta]_{A\oplus
B}=[\rho_B(\xi),\rho_B(\eta)]_\dev.$$ By Lemma \ref{lem:rho} and
Conditions (2) and (3) of Definition \ref{Def:Ebialgebroid}, we have
$$
[\rho_A (X),\rho_B (\xi)]_\dev =\rho_B (\Liederivative_X \xi)-\rho_A
(\Liederivative_\xi X)+\rho_A \dM^B\langle X,\xi\rangle_E.
$$
We also have, \begin{eqnarray*} (\rho_A+ \rho_B)[X, \xi]_{A\oplus
B}&=&(\rho_A+ \rho_B)(\Liederivative_X\xi-\Liederivative_\xi
X+\dM^B\langle X,\xi\rangle_E)\\&=&\rho_B (\Liederivative_X
\xi)-\rho_A (\Liederivative_\xi X)+\rho_A \dM^B\langle
X,\xi\rangle_E.
 \end{eqnarray*}
So we  get Property  (EC-1). Since $(\rho_A+
\rho_B)^\star=\rho_A^\star+\rho_B^\star$ and $\rho_A^\star(\jetd
u)=\dM^Au$, we have
$$
[X+\xi,X+\xi]_{A\oplus B}=(\dA+\dB)\langle X,\xi\rangle_E=(\rho_A+
\rho_B)^\star\circ\jetd\langle X,\xi\rangle_E\,,
$$
and  Property (EC-2) follows. Property (EC-3) is straightforward.
Property (EC-4) follows from the fact that $\rho_A^\star\jet
E\subset B$ and $\rho_B^\star\jet E\subset A$. Property (EC-5)
follows from Condition (3) of Definition \ref{Def:Ebialgebroid}.
This proves that  $(A\oplus
B,\ppairingE{\cdot,\cdot},[\cdot,\cdot]_{A\oplus B},\rho_A+\rho_B)$
is an $E$-Courant algebroid. \qed \vspace{3mm}

Conversely, for an $E$-Courant algebroid $\huaK$, suppose that there
are two transverse Dirac structures $A$ and $B$ in $\huaK$ such that
$\huaK=A\oplus B$, we want to show that $(A,B)$ is an $E$-Lie
bialgebroid.

Evidently, $A$ and $B$ are mutually $E$-dual vector bundles, with
the pairing $$\langle
X,\xi\rangle_E=2(X,\xi)_E,~X\in\Gamma(A),~\xi\in\Gamma(B) .$$ By
Proposition \ref{pro:rep}, both $A$ and $B$ are Lie algebroids whose
anchors are, respectively, $a=\jd\circ\rho\mid_A$ and
$a_B=\jd\circ\rho\mid_B$. In the meantime, there are representations
 $\rho_A=\rho\mid_A$ of $A$ and $\rho_B=\rho\mid_B$ of $B$ on $E$.
 The associated coboundary operators $\dA:\Omega^\bullet( A, E)\lon
\Omega^{\bullet+1}( A, E)$ and $\dB:\Omega^\bullet( B, E)\lon
\Omega^{\bullet+1}( B, E)$ are standard.

 Similar to the result
in \cite{LWXmani} for Courant algebroids, we have the $E$-Courant
algebroid analogue:
\begin{thm}\label{thm:manin}
If an $E$-Courant algebroid has a  decomposition  $\huaK=A\oplus B$,
where $A$ and $B$ are transverse Dirac structures, then $(A,B)$ is
an $E$-Lie bialgebroid.
\end{thm}
\pf By Property (EC-2), for any $X\in \Gamma(A),~\xi\in \Gamma(B)$,
we have
\begin{equation}\label{eqn skew}
{[X,\xi]}_{\huaK}=-[\xi,X]_{\huaK}+(\dA+\dB)\langle X,\xi\rangle_E.
\end{equation}
By (\ref{eqn skew}) and Property (EC-3), we have
\begin{eqnarray*}
[X,\xi]_{\huaK}&=&\Liederivative_X\xi-\Liederivative_{\xi}X+\dB\langle
X,\xi\rangle_E, \\
{[\xi,X]}_{\huaK}
&=&\Liederivative_{\xi}X-\Liederivative_X\xi+\dA\langle
X,\xi\rangle_E.
\end{eqnarray*}
So  we have $\Liederivative_X\xi\in \Gamma(B)$ and
$\Liederivative_\xi X\in \Gamma(A)$, and one gets
$\dA(\Gamma(B))\subset\Gamma(\Hom(\wedge^2 A,E)_B)$, ~
$\dB(\Gamma(A))\subset\Gamma(\Hom(\wedge^2B,E)_A)$. Furthermore, we
have $\rhowx\circ\jetd=(\rhoAWuXing+\rhoBWuXing)\circ\jetd=\dA+\dB$.
By Property (EC-4), we have $\rhowx(\jet E)\subset\huaK$ and hence
$\dA(\Gamma(E))\subset \Gamma(B)$ and $\dB (\Gamma(E))\subset
\Gamma(A)$. Therefore, $((A,\rho_A);(B,\rho_B))$ is an $E$-dual
pair.

Under the decomposition $\huaK=A\oplus B$, for sections
$e_i=X_i+\xi_i\in \Gamma(\huaK)$, $i=1,2$, the bracket
$[e_1,e_2]_{\huaK}$ is given by (\ref{bracket1}). By Lemma
\ref{lem:jacobi} and the fact that $ \Gamma(\huaK)$ is a Leibniz
algebra, we have   $J_1+c.p.+J_2=0 $. Moreover, we have
\begin{equation}\label{eqn:a d}
[\rho_B(\xi),\rho_A(X)]_\dev=\rho[\xi,X]_{\huaK}=\rho_A(\Liederivative_{\xi}X)-\rho_B(\Liederivative_X\xi)+
\rho_B\circ\dA\langle\xi,X\rangle_E.
\end{equation}
By  Lemma \ref{lem:rho}, we have
\begin{equation}\label{eqn:condition2}
\langle\Liederivative_{\dB u}\xi+[\dA u,\xi],X\rangle_E=0.
\end{equation}
The nondegeneracy of the $E$-valued pairing implies  that
$\Liederivative_{\dB u}\xi=-\Liederivative_{\dA u}\xi$, i.e.
 Condition (2) of Definition
\ref{Def:Ebialgebroid}.

Finally by  (\ref{eqn:condition2}), we have  $J_2=0$, which implies
that $J_1+c.p.=0$. In particular, if we take $e_1=X_1,~e_2=X_2$ and
$e_3=\xi_3$, we have $\dB[X_1,X_2] -\Liederivative_{X_1}\dB
X_2+\Liederivative_{X_2}\dB X_1=0$, which is equivalent to Condition
(1) of Definition \ref{Def:Ebialgebroid}.

By Property  (EC-5), we have
$\rho_A\circ\rhoBWuXing=-\rho_B\circ\rhoAWuXing$. So Condition (3)
of Definition \ref{Def:Ebialgebroid} holds. In summary, $(A,B)$ is
an $E$-Lie bialgebroid. \qed \vspace{3mm}

Finally, we give some   examples of $E$-Lie bialgebroids.

$\bullet$ {\bf{The $T^*M$-Lie bialgebroid $(\jet A, \jet (A^*))$}}

 For any Lie bialgebroid $(A,A^*)$, there  associates a Courant
algebroid structure on $A\oplus A^*$. By Theorem \ref{thm:jet C},
$\jet (A\oplus A^*)=\jet A\oplus \jet (A^*)$ is a $T^*M$-Courant
algebroid. It is easily seen that both $\jet A$ and $ \jet (A^*)$
are transverse Dirac structures. Thus by Theorem \ref{thm:manin}, we
have
\begin{pro}
For any Lie bialgebroid $(A,A^*)$, $(\jet A, \jet (A^*))$ is a
$T^*M$-Lie bialgebroid.
\end{pro}

$\bullet$ {\bf{The $E$-Lie bialgebroid $(\dev E,\jet E)$ induced by
a Lie algebroid $(E,[\cdot,\cdot]_E,a)$}}

Recall Theorem \ref{Thm:ElocalLie} where it is shown that for any
$E$-Lie bialgebroid, there is a
  Lie algebroid structure (local Lie algebra structure) on $E$ if
  $\mathrm{rank}E\geq2$ ($\mathrm{rank}E=1$). There is also  a canonical $E$-Lie bialgebroid $(\dev
E,\jet E)$ coming from a given Lie algebroid
$(E,[\cdot,\cdot]_E,a)$. If $\mathrm{rank}E=1$, similar conclusion
holds for any
 local Lie algebra $(E,[\cdot,\cdot]_E)$.

In fact,  if
 $(E,[\cdot,\cdot]_E,a)$ is a  Lie
algebroid (or  $(E,[\cdot,\cdot]_E)$ is a local Lie algebra when
$\mathrm{rank}E=1$), we can define a skew-symmetric bundle map
$\pi:\jet E\longrightarrow\dev E$ by setting
$$
\pi(\jetd u)(v)=[u,v]_E,\quad \forall~u,~v\in\Gamma(E).
$$(This appeared in \cite{clomni}.)
Moreover, we introduce a bracket $[\cdot,\cdot]_\pi$ on $\Gamma(\jet
E)$ by setting
\begin{equation}
\pibracket{\mu,\nu}\defbe
\Liederivative_{\pi(\mu)}\nu-\Liederivative_{\pi(\nu)}\mu-
\jetd\conpairing{\pi(\mu),\nu}_E.
\end{equation}
One is able to check that   the graph of $\pi$ is a Dirac structure
and
  $(\jet E,[\cdot,\cdot]_\pi,\jd\circ\pi)$ is a Lie
algebroid equipped with a representation $\pi$. Therefore, we obtain
an   $E$-Lie bialgebroid structure on $(\dev E,\jet E)$.

For a Poisson manifold $(M,\pi)$, there is a canonical Lie
bialgebroid $(TM, T^*M)$. For a Jacobi manifold $(M,X,\Lambda)$,
there is a canonical generalized Lie bialgebroid $(TM\oplus\mathbb
R,T^*M\oplus\mathbb R)$. Similarly,  we have
\begin{pro}
For any  Lie algebroid $(E,[\cdot,\cdot]_E,a)$, there is a canonical
$E$-Lie bialgebroid $(\dev E,\jet E)$. If $\mathrm{rank}E=1$, the
conclusion holds for any
 local Lie algebra structure $(E,[\cdot,\cdot]_E)$.
\end{pro}

 $\bullet$
{\bf{Dirac structures}}

For an $E$-Lie bialgebroid $(A,B)$, suppose that
$H\in\Gamma(\Hom(\wedge^2B,E)_A)$. Treated as  a bundle map
$H:B\longrightarrow A $, $H$ has its graph
$\huaG_H=\{H\xi+\xi|~\forall~\xi\in B\}\subset A\oplus B$.

\begin{thm}
The graph $\huaG_H$ is a Dirac structure if and only if $H$
satisfies the following Maurer-Cartan type equation:
$$
\dB H+\frac{1}{2}[H,H]=0.
$$
\end{thm}
\pf The property that $\huaG_H$ is isotropic is equivalent to the
condition  that $H\in\Gamma(\Hom(\wedge^2B,E)_A)$. By
(\ref{bracket1}), we have
$$
[H\xi,\eta]+[\xi,H\eta]=[\xi,\eta]_H+\Liederivative_\xi
H\eta-\Liederivative_\eta H\xi+\dB (H(\xi,\eta)).
$$
By (\ref{eqn: H}), we have
$$
[H\xi+\xi,H\eta+\eta]=\Liederivative_\xi H\eta-\Liederivative_\eta
H\xi+\dB
(H(\xi,\eta))+H[\xi,\eta]_H+\frac{1}{2}[H,H](\xi,\eta)+[\xi,\eta]+[\xi,\eta]_H.
$$
So $\Gamma(\huaG_H)$ is closed under the bracket $[~,~]$  if and
only if for any $\xi,~\eta\in \Gamma(B)$,
\begin{equation}\label{eqn:H closed}
H[\xi,\eta]=\Liederivative_\xi H\eta-\Liederivative_\eta H\xi+\dB
(H(\xi,\eta))+\frac{1}{2}[H,H](\xi,\eta).
\end{equation}
On the other hand,
\begin{eqnarray*}
(\dB
H)(\xi,\eta,\vartheta)&=&\rho_B(\xi)\conpairing{H\eta,\vartheta}_E-\rho_B(\eta)\conpairing{H\xi,\vartheta}_E
+\rho_B(\vartheta)\conpairing{H\xi,\eta}_E\\&&+\conpairing{H\vartheta,[\xi,\eta]}_E-\conpairing{H\eta,[\xi,\vartheta]}_E
+\conpairing{H\xi,[\eta,\vartheta]}_E\\
&=&\conpairing{\Liederivative_\xi H\eta-\Liederivative_\eta
H\xi-H[\xi,\eta]+\dB (H(\xi,\eta)),\vartheta}_E.
\end{eqnarray*}
Therefore,  (\ref{eqn:H closed}) is equivalent to the condition that
$$ (\dB H)(\xi,\eta)+\frac{1}{2}[H,H](\xi,\eta)=0,
$$
or equivalently,
$$
\dB H+\frac{1}{2}[H,H]=0. \qed
$$

In particular, if $\dB=0$ (i.e. $B$ is a trivial Lie algebroid), the
graph $\huaG_H$ is a Dirac structure if and only if $[H,H]=0$, or
$H[\xi,\eta]_H=[H(\xi),H(\eta)]$.

\end{document}